# Goodness-of-fit tests for Markovian time series models: Central limit theory and bootstrap approximations

MICHAEL H. NEUMANN[1] and EFSTATHIOS PAPARODITIS[2]

[1]*Friedrich Schiller Universität Jena, Institut für Stochastik, Ernst Abbe Platz 2, D–07743 Jena, Germany. E-mail: mneumann@mathematik.uni-jena.de*
[2]*University of Cyprus, Department of Mathematics and Statistics, P.O. Box 20537, CY–1678 Nicosia, Cyprus. E-mail: stathisp@ucy.ac.cy*

New goodness-of-fit tests for Markovian models in time series analysis are developed which are based on the difference between a fully nonparametric estimate of the one-step transition distribution function of the observed process and that of the model class postulated under the null hypothesis. The model specification under the null allows for Markovian models, the transition mechanisms of which depend on an unknown vector of parameters and an unspecified distribution of i.i.d. innovations. Asymptotic properties of the test statistic are derived and the critical values of the test are found using appropriate bootstrap schemes. General properties of the bootstrap for Markovian processes are derived. A new central limit theorem for triangular arrays of weakly dependent random variables is obtained. For the proof of stochastic equicontinuity of multidimensional empirical processes, we use a simple approach based on an anisotropic tiling of the space. The finite-sample behavior of the proposed test is illustrated by some numerical examples and a real-data application is given.

## 1. Introduction

The analysis of time series is often based on parametric or semi-parametric model assumptions that must be tested in practice. An important class of stochastic processes used in modelling time series is that of Markov type, which are described by a specification of the transition kernel, usually involving a finite-dimensional parameter vector and possibly a partial specification of the distribution of some innovations.

Many available tests of parametric or semiparametric models in time series analysis originate from corresponding frameworks with independent and identically distributed data and are based on the difference between the stationary distribution or the regression/autoregression function and corresponding empirical counterparts. Tests proposed by Bierens (1982) and McKeague and Zhang (1994) focus on the conditional mean function. For the autoregression function, an approach leading to distribution-free tests using a martingale transformation has been proposed by Koul and Stute (1999). For an overview, see also Delgado and González Manteiga (2001). Sometimes, the autoregression function that corresponds to a hypothetical model is not explicitly known, which causes problems for testing schemes based on it. More importantly, however, a model check based on the stationary distribution or the autoregression function is not able to detect all types of departures from a hypothetical model. It might happen that the autoregression functions and/or the conditional variance functions of two processes are similar or identical, whereas their conditional distributions are essentially different.

In the present paper, we derive tests of the validity of a Markov model by directly comparing the hypothetical, model-based conditional distribution with its model-free estimated empirical counterpart. Let $\mathbf{X} = (X_t)_{t \in \mathbb{Z}}$ be the process considered and denote by $P_\mathbf{X}$ the law of $\mathbf{X}$. Denote by $\mathcal{M}$ the class of Markov processes of order less than or equal to $p$, that is,

$$\mathcal{M} = \{P_\mathbf{X} : P(X_t \in B | \sigma(X_s, s < t)) = P(X_t \in B \mid \mathbb{X}_{t-1}) \ \forall B \in \mathcal{B}, \forall t \in \mathbb{Z}\},$$

where $\mathbb{X}_{t-1} = (X_{t-1}, \ldots, X_{t-p})'$. The problem we consider in this paper is that of testing the hypothesis

$$H_0 : P_\mathbf{X} \in \mathcal{M}_0 \cap \mathcal{M}$$

against the alternative

$$H_1 : P_\mathbf{X} \in \mathcal{M} \setminus \mathcal{M}_0,$$

where $\mathcal{M}_0 \subset \mathcal{M}$ denotes the class of $p$th order Markov processes described by

$$\mathcal{M}_0 = \{P_\mathbf{X} : X_t = G(\mathbb{X}_{t-1}, \theta, \varepsilon_t), \varepsilon_t \sim F_\varepsilon \text{ i.i.d.}; \theta \in \Theta, F_\varepsilon \in \mathcal{F}_\varepsilon\}.$$

In the above notation, $G : \mathbb{R}^p \times \Theta \times \mathbb{R} \to \mathbb{R}$ is some known function depending on an unknown parameter vector $\theta$ and $(\varepsilon_t)_{t \in \mathbb{Z}}$ is a sequence of i.i.d. innovations with distribution $F_\varepsilon$ belonging to some appropriate class of distribution functions denoted by $\mathcal{F}_\varepsilon$. Note that, under the null hypothesis, the transition kernel generating $X_t$ is known up to the finite-dimensional parameter vector $\theta$ and the distribution function $F_\varepsilon$ of the innovations. Testing problems fitting into this framework are discussed in Section 2. In this paper, an appropriate test statistic for the above testing problem is obtained which is based on the supremum deviation of a fully nonparametric estimator of the conditional distribution function $F_{X_t | \mathbb{X}_{t-1}}$ from its model-based one. If $H_0$ were a simple hypothesis, that is, if $\theta$ and $F_\varepsilon$ were specified under $H_0$, then we could consider as a starting point for a test the deviation process given by

$$U_n^{(0)}(x, y) = \frac{1}{\sqrt{n}} \sum_{t=1}^n I(\mathbb{X}_{t-1} \preceq x)[I(X_t \leq y) - F_{H_0}(y \mid \mathbb{X}_{t-1})],$$

$$(x,y) \in \mathbb{R}^p \times \mathbb{R}, \tag{1.1}$$

where $F_{H_0}(y \,|\, \mathbb{X}_{t-1})$ denotes the one-step transition distribution function of the Markov model class postulated under $H_0$ and $x \preceq y$ means $x_i \leq y_i$ for all $i = 1, \ldots, p$. Note that $U_n^{(0)}(x,y)$ measures, for every point $(x,y)$, the difference between the sample one-step transition distribution function and its model-based version under the null. As we will see in the sequel, appropriate specifications of the function $F_{H_0}(y \,|\, \mathbb{X}_{t-1})$ under $H_0$ lead to useful test statistics. We consider as important cases $\mathrm{AR}(p)$ or $\mathrm{ARCH}(p)$ processes and show that the above deviation process can be asymptotically approximated by a Gaussian process. Note that, although $\mathcal{M}_0$ is general enough, it may be of interest in some applications to test more restricted versions of $\mathcal{M}_0$. For instance, in ARCH modelling, it is not uncommon to impose specific parametric assumptions on the distribution $F_\varepsilon$ of the innovations (cf. Engle (1982) and Bollerslev (1987)). Our inference procedure also allows for the testing of hypotheses of this type; see the discussion at the end of Section 3.

Since the parameters of the limiting Gaussian process of the statistic considered depend on the actual stochastic process in a complicated way, we approximate the null distribution of the test statistic by means of a model-based bootstrap approach. For the cases considered, the most natural bootstrap scheme is a model-based bootstrap without an additional smoothing of the residuals, that is, a bootstrap based on i.i.d. residuals with a discrete distribution. It is known that such processes may fail to satisfy classical mixing conditions; see, for example, Rosenblatt (1980). However, using the alternative concept of weak dependence introduced by Doukhan and Louhichi (1999), we are able to establish some general properties of such bootstrap schemes which may be of interest in their own right. In a different context and for inference problems different from those considered here, bootstrap methods related to that proposed in this paper have been considered by, among others, Basawa, Green, McCormick and Taylor (1990), Rajarshi (1990), Paparoditis and Politis (2002) and Bühlmann (2002).

Testing problems based on the conditional distribution function have attracted some interest in recent years. For i.i.d. observations, specification tests of a parametric hypothesis concerning the conditional distribution have been considered by Andrews (1997) and Stinchcombe and White (1998). For dependent, mostly-mixing observations, tests based on the conditional distribution function have also been investigated by some authors. Li and Tkacz (2001) proposed a test based on an $L^2$-type distance between the nonparametrically estimated conditional density and its model-based parametric counterpart. Corradi and Swanson (2001) and Bai (2003) considered a Kolmogorov-type test. Testing problems associated with linear restrictions on the conditional distribution function have been considered by Inoue (1999). Although some of the above approaches deal with testing problems similar to our own, none of these papers covers the more general case described by our null hypothesis $H_0$. Furthermore, our analysis is based on the alternative concept of weak dependence, which appears to be quite important in the current context. We investigate basic properties of model-based bootstrap approaches and, as a technical prerequisite for our analysis, we prove a central limit theorem for triangular arrays of weakly dependent random variables which does not require moment conditions beyond Lindeberg's.

The paper is organized as follows. In Section 2, we precisely state our assumptions on the underlying stochastic process and discuss some interesting examples of Markovian models which fit into our testing framework. In Section 3, the test statistic used is presented and its asymptotic behavior under validity of the null hypothesis is established. Bootstrap approximations to the distribution of the test statistic under the null are investigated in Section 4. Some numerical examples illustrating our theoretical analysis, as well as an application to financial data, are given in Section 5. In Section 6, we prove a central limit theorem for triangular schemes of weakly dependent random variables. Section 7 contains the proofs of all statistical results and a lemma which provides the major step in a proof of stochastic equicontinuity.

## 2. Assumptions and examples of processes

We assume that observations $X_{1-p}, \ldots, X_n$ from a real-valued, stationary process $\mathbf{X} = (X_t)_{t \in \mathbb{Z}}$ are available. The null hypothesis is that $\mathbf{X}$ is a Markov process of order $p$ with a particular form of the conditional distribution, that is, $P_\mathbf{X}$ belongs to $\mathcal{M}_0 \cap \mathcal{M}$. Recall that $G: \mathbb{R}^p \times \Theta \times \mathbb{R} \to \mathbb{R}$ is some known function depending on an unknown parameter vector $\theta$ and $(\varepsilon_t)_{t \in \mathbb{Z}}$ is a sequence of i.i.d. innovations with distribution $F_\varepsilon$. Since our test statistic below requires estimates of $\theta$ and $F_\varepsilon$, we must be more specific about the models to be considered and it turns out that asymptotic properties must be derived in a case-by-case manner. We will focus our attention on the following classes of processes.

### 2.1. AR($p$) processes

Here, we suppose that, under $H_0$, the following condition is satisfied.

(A1) The process $\mathbf{X} = (X_t)_{t \in \mathbb{Z}}$ obeys the model equation

$$X_t = \theta' \mathbb{X}_{t-1} + \varepsilon_t, \tag{2.1}$$

where $(\varepsilon_t)_{t \in \mathbb{Z}}$ is a sequence of i.i.d. innovations with $\mathrm{E}\varepsilon_t = 0$, $0 < \mathrm{E}\varepsilon_t^2 = \sigma^2$, $\mathrm{E}\varepsilon_t^4 < \infty$ and $\theta = (\theta_1, \theta_2, \ldots, \theta_p)' \in \Theta$ and

$$\Theta = \{\theta \in \mathbb{R}^p : 1 - \theta_1 z - \cdots - \theta_p z^p \neq 0 \text{ for all } z \in \mathbb{C} \text{ with } |z| \leq 1\}.$$

In this case, the function $G$ is given as $G(\mathbb{X}_{t-1}, \theta, \varepsilon_t) = \theta' \mathbb{X}_{t-1} + \varepsilon_t$, while $P(X_t \leq y \mid \mathbb{X}_{t-1}) = F_\varepsilon(w(\mathbb{X}_{t-1}, y, \theta))$ with $w(\mathbb{X}_{t-1}, y, \theta) = y - \theta' \mathbb{X}_{t-1}$.

It is well known that (2.1) has, under (A1), a unique stationary solution which has a representation as a causal linear (MA($\infty$)-) process, that is, $X_t = \sum_{k=0}^\infty \alpha_k \varepsilon_{t-k}$, where $|\alpha_k| \leq K \rho^k$ for some $K < \infty$, $\rho < 1$; see, for example, Brockwell and Davis (1991), page 85. Adapting the proof of Lemma 9 in Doukhan and Louhichi (1999), it can be shown that (A1) implies, for $s_1 < \cdots < s_u < t_1 < \cdots < t_v$ and arbitrary measurable functions $g: \mathbb{R}^u \longrightarrow \mathbb{R}$, $h: \mathbb{R}^v \longrightarrow \mathbb{R}$ with $\mathrm{E}g^2(X_{s_1}, \ldots, X_{s_u}) < \infty$, $\mathrm{E}h^2(X_{t_1}, \ldots, X_{t_v}) < \infty$,

$$|\mathrm{cov}(g(X_{s_1}, \ldots, X_{s_u}), h(X_{t_1}, \ldots, X_{t_v}))|$$

$$\leq \sqrt{\mathrm{E}g(X_{s_1},\ldots,X_{s_u})^2}\,\mathrm{Lip}(h)\sqrt{\mathrm{E}\varepsilon_0^2}\sum_{k=t_1-s_u}^{\infty}|\alpha_k|[(k-(t_1-s_u)+1)\wedge v]. \tag{2.2}$$

Furthermore, it can also be shown that, for $s_1 < \cdots < s_u < t_1 \leq t_2$, $1 \leq k_1, k_2 \leq p$, the following inequalities hold true: for any measurable $g\colon\mathbb{R}^u \longrightarrow \mathbb{R}$ with $\mathrm{E}g^2(X_{s_1},\ldots,X_{s_u}) < \infty$,

$$|\mathrm{cov}(g(X_{s_1},\ldots,X_{s_u}),X_{t_1-k_1}\varepsilon_{t_1})| \leq C\sqrt{\mathrm{E}g^2(X_{s_1},\ldots,X_{s_u})}\rho_{t_1-k_1-s_u} \tag{2.3}$$

and, for any measurable $g\colon\mathbb{R}^u \longrightarrow \mathbb{R}$ with $\|g\|_\infty = \sup_{x\in\mathbb{R}^u}|g(x)| < \infty$,

$$|\mathrm{cov}(g(X_{s_1},\ldots,X_{s_u}),X_{t_1-k_1}X_{t_2-k_2}\varepsilon_{t_1}\varepsilon_{t_2})| \leq C\|g\|_\infty\rho_{\min\{t_1-k_1-s_u,t_2-k_2-s_u\}}, \tag{2.4}$$

where $\rho_r = \sqrt{\sum_{k=r}^\infty |\alpha_k|}$, that is, a weak dependence condition similar to those in Doukhan and Louhichi (1999) is fulfilled.

## 2.2. ARCH($p$) processes

The class of autoregressive, conditionally heteroscedastic (ARCH) processes was introduced by Engle (1982). In this case, our null hypothesis means that the process fulfills the following condition.

(A1′) The process $\mathbf{X} = (X_t)_{t\in\mathbb{Z}}$ is stationary and obeys the model equation

$$X_t = \sqrt{\theta_0 + \theta_1 X_{t-1}^2 + \cdots + \theta_p X_{t-p}^2}\,\varepsilon_t, \tag{2.5}$$

where $\theta = (\theta_0,\theta_1,\ldots,\theta_p)' \in \Theta$ and $\Theta = \{\theta \in \mathbb{R}^{p+1}: \theta_0 > 0, \theta_i \geq 0, i = 1,\ldots,p$ and $\sum_{i=1}^p \theta_i < 1\}$. Furthermore, $(\varepsilon_t)_{t\in\mathbb{Z}}$ is a sequence of i.i.d. innovations with $\mathrm{E}\varepsilon_t = 0$, $\mathrm{E}\varepsilon_t^2 = 1$ and $\mathrm{E}\varepsilon_t^8 < \infty$.

In this case, the function $G$ is given as $G(\mathbb{X}_{t-1},\theta,\varepsilon_t) = \sqrt{\theta_0 + \theta_1 X_{t-1}^2 + \cdots + \theta_p X_{t-p}^2}\,\varepsilon_t$. Moreover, we have that $P(X_t \leq y \mid \mathbb{X}_{t-1}) = F_\varepsilon(w(\mathbb{X}_{t-1},y,\theta))$ with $w(\mathbb{X}_t,y,\theta) = y/\sqrt{\theta_0 + \theta_1 X_{t-1}^2 + \cdots + \theta_p X_{t-p}^2}$. Milhøj (1985) obtained a representation of the unique stationary solution to (2.5) as $\widetilde{X}_t = \varepsilon_t \cdot (\theta_0 \cdot \sum_{k=0}^\infty M(t,k))^{1/2}$, where $M(t,0) = 1$, $M(t,k) = \sum_{a_1,\ldots,a_k=1}^p \prod_{i=1}^k \theta_{a_i}\varepsilon_{t-a_1-\cdots-a_i}^2$. (Obviously, $\sum_{i=1}^k$ in the second display after equation (2.1) in Milhøj (1985) should read $\prod_{i=1}^k$.) The process $(\widetilde{X}_t)_{t\in\mathbb{Z}}$ is both weakly and strictly stationary. To deal with this process, we will exploit the following weak dependence property.

**Lemma 2.1.** *Suppose that* (A1′) *is fulfilled. There then exist some $\rho < 1$, $C < \infty$, such that, for all $s_1 < \cdots < s_u < t_1 < \cdots < t_v$ and arbitrary measurable functions $g\colon\mathbb{R}^u \longrightarrow \mathbb{R}$,* $h : \mathbb{R}^v \longrightarrow \mathbb{R}$ *with* $\mathrm{E}g^2(X_{s_1},\ldots,X_{s_u}) < \infty$ *and* $\mathrm{E}h^2(X_{t_1},\ldots,X_{t_v}) < \infty$,

$$\begin{aligned} &|\mathrm{cov}(g(X_{s_1},\ldots,X_{s_u}), h(X_{t_1},\ldots,X_{t_v}))| \\ &\quad \leq C\sqrt{\mathrm{E}g^2(X_{s_1},\ldots,X_{s_u})}\mathrm{Lip}(h)\rho^{t_1-s_u}. \end{aligned} \quad (2.6)$$

*Furthermore, for* $s_1 < \cdots < s_u < t_1 \leq t_2$, $1 \leq k_1, k_2, k_3, k_4 \leq p$, *the following inequalities hold true: for any measurable* $g : \mathbb{R}^u \longrightarrow \mathbb{R}$ *with* $\mathrm{E}g^2(X_{s_1},\ldots,X_{s_u}) < \infty$,

$$\begin{aligned} &|\mathrm{cov}(g(X_{s_1},\ldots,X_{s_u}), X_{t_1-k_1}^2 X_{t_1-k_2}^2(\varepsilon_{t_1}^2 - 1))| \\ &\quad \leq C\sqrt{\mathrm{E}g^2(X_{s_1},\ldots,X_{s_u})}\rho^{t_1-s_u} \end{aligned} \quad (2.7)$$

*and, for any measurable* $g : \mathbb{R}^u \longrightarrow \mathbb{R}$ *with* $\|g\|_\infty = \sup_{x \in \mathbb{R}^u} |g(x)| < \infty$,

$$\begin{aligned} &|\mathrm{cov}(g(X_{s_1},\ldots,X_{s_u}), X_{t_1-k_1}^2 X_{t_1-k_2}^2 X_{t_2-k_3}^2 X_{t_2-k_4}^2(\varepsilon_{t_1}^2 - 1)(\varepsilon_{t_2}^2 - 1))| \\ &\quad \leq C\|g\|_\infty \rho^{t_1-s_u}. \end{aligned} \quad (2.8)$$

## 2.3. Markov processes driven by diffusions

At this point, we would like to explicitly mention two particular classes of processes which are of interest in financial mathematics.

### 2.3.1. Vasicek model

Merton (1971) proposed to model interest rate processes by an Ornstein–Uhlenbeck process (in continuous time),

$$\mathrm{d}X_t = \theta_1(\theta_2 - X_t)\,\mathrm{d}t + \theta_3\,\mathrm{d}W_t, \quad (2.9)$$

where $\theta_1, \theta_2, \theta_3 > 0$ and $(W_t)_{t \in \mathbb{Z}}$ is a standard Wiener process. This model was further considered by Vasicek (1977). Assume that we observe this process at equidistant design points $\Delta t$, where $t = 0, 1, \ldots, n$. These observations then form a Markov process with stationary and conditionally normal distributed increments where, for $s < t$, $\mathrm{E}(X_t \mid X_s) = \theta_2 + (X_s - \theta_2)\mathrm{e}^{-\theta_1(t-s)}$ and $\mathrm{var}(X_t \mid X_s) = \theta_3^2(1 - \exp\{-2\theta_1(t-s)\})/(2\theta_1)$. Hence, introducing appropriate innovations $\varepsilon_t$, we can rewrite $(X_{\Delta t})_{t \in \mathbb{Z}}$ in the form $X_{\Delta t} = G(X_{\Delta(t-1)}, \theta, \varepsilon_t)$, where $G(X_{\Delta(t-1)}, \theta, \varepsilon_t) = [\theta_2 + (X_{\Delta(t-1)} - \theta_2)\mathrm{e}^{-\theta_1 \Delta}] + \varepsilon_t$ and $(\varepsilon_t)_{t \in \mathbb{Z}}$ is a sequence of independent normally distributed variables with zero mean and variance $\theta_3^2(1 - \exp\{-2\theta_1(t-s)\})/(2\theta_1)$. Accordingly, the process $(Y_t)_{t \in \mathbb{Z}}$ with $Y_t = X_{\Delta t} - \theta_2$ is an AR(1) process with parameter $\mathrm{e}^{\theta_1 \Delta}$ which satisfies conditions analogous to (2.2) to (2.4) above.

*2.3.2. Cox–Ingersoll–Ross model*

Also for the purpose of modelling interest rates, Cox, Ingersoll and Ross (1985) proposed the specification

$$\mathrm{d}X_t = \theta_1(\theta_2 - X_t)\,\mathrm{d}t + \theta_3\sqrt{X_t}\,\mathrm{d}W_t, \tag{2.10}$$

where $\theta_1, \theta_2, \theta_3 > 0$ and $(W_t)_{t\geq 0}$ is a standard Wiener process. Again, the values of the process at equidistant time points form a stationary Markovian process. An explicit description of the conditional distribution of $X_{\Delta t}$ given $X_{\Delta(t-1)}$ can be found in Cox *et al.* (1985), page 391. In the special case where $q = 4\theta_1\theta_2/\theta_3^2$ is an integer, it follows that the conditional distribution of $Y_t = cX_{\Delta t}$ $[c = 4\theta_1/(\theta_3(1-\mathrm{e}^{-\theta_1\Delta}))]$ given $Y_{t-1}$ is noncentral chi-square with $q$ degrees of freedom and parameter of noncentrality $Y_{t-1}\mathrm{e}^{-\theta_1\Delta}$.

Let $G(\lambda, x) = F^{-1}_{\chi^2_q(\exp\{-\theta_1\Delta\}\lambda)}(x)$. With an appropriate sequence $(\varepsilon_t)_{t\in\mathbb{Z}}$ of independent uniform$(0,1)$-distributed innovations, we can write $(Y_t)_{t\in\mathbb{Z}}$ as $Y_t = G(Y_{t-1}, \theta, \varepsilon_t)$, $\theta = (\theta_1, \theta_2, \theta_3)'$. From $\mathrm{E}|G(0, \theta, \varepsilon_t)| < \infty$ and $\mathrm{E}|G(u, \theta, \varepsilon_t) - G(v, \theta, \varepsilon_t)| = |\mathrm{E}G(u, \theta, \varepsilon_t) - \mathrm{E}G(v, \theta, \varepsilon_t)| \leq \exp\{-\theta_1\Delta\}|u-v|$, it follows, analogously to Lemma 10 of Doukhan and Louhichi (1999), that conditions such as (2.2) to (2.4) are fulfilled.

## 3. The test statistics and its limit distribution

Consider testing the hypothesis $H_0$ of interest. Since $P(X_t \leq y \mid \mathbb{X}_{t-1} = x)$ cannot be consistently estimated in general, we construct our test statistic from *cumulative* versions of the hypothetical transition probabilities and model-free estimators thereof. As mentioned in the Introduction, the basic idea needed to construct an appropriate deviation process is to consider the difference between a fully nonparametric version of the one-step transition distribution function and its parametric version postulated under $H_0$, that is, to consider the basic deviation process $U_n^{(0)}(x, y)$ given in (1.1). To specify $F_{H_0}(y \mid \mathbb{X}_{t-1})$ given in this equation, that is, to specify the one-step transition distribution function under the null hypothesis, we proceed as follows. Since the null hypothesis is only partially specified, that is, $\theta$ and $F_\varepsilon$ are unknown, we replace these unknown quantities by their corresponding sample estimates. In particular, and in order to deal with the uncertainty introduced by the fact that $\theta$ is unknown, we assume the following

(A2) The sequence of estimators $\widehat{\theta}$ admits the expansion

$$\widehat{\theta} - \theta = \frac{1}{n}\sum_{t=1}^{n} l(\mathbb{X}_{t-1}, X_t; \theta) + o_P(n^{-1/2}),$$

where $l(\cdot, \cdot; \cdot)$ is a measurable function from $\mathbb{R}^p \times \mathbb{R} \times \Theta$ to $\mathbb{R}^k$ with $\mathrm{E}_\theta l(\mathbb{X}_{t-1}, X_t; \theta) = (0, \ldots, 0)'$ and $\mathrm{E}_\theta \|l(\mathbb{X}_{t-1}, X_t; \theta)\|^2 < \infty$.

Note that in the AR$(p)$ case, assumption (A2) is satisfied with $k = p$, $l(\mathbb{X}_{t-1}, X_t; \theta) = \Gamma_p^{-1}\mathbb{X}_{t-1}\varepsilon_t$ and $\Gamma_p = \Gamma_p(\theta) = \mathrm{E}_\theta(\mathbb{X}_{t-1}\mathbb{X}'_{t-1})$ if $\widehat{\theta}$ is the commonly used least-squares or Yule–Walker estimator (cf. Brockwell and Davis (1991)). For the linear ARCH$(p)$ case, (A2) is, for instance, satisfied if $\widehat{\theta}$ is the least-squares estimator of $\theta = (\theta_0, \theta_1, \ldots, \theta_p)'$. In this case, $k = p+1$, $l(\mathbb{X}_{t-1}, X_t; \theta) = C_p^{-1} \mathbb{Y}_{t-1} \mathbb{Y}'_{t-1} \theta(\varepsilon_t^2 - 1)$, where $C_p = C_p(\theta) = \mathrm{E}_\theta(\mathbb{Y}_{t-1} \mathbb{Y}'_{t-1})$ and $\mathbb{Y}_{t-1} = (1, X_{t-1}^2, \ldots, X_{t-p}^2)'$.

To deal with the fact that $F_\varepsilon$ is unknown, we replace $F_\varepsilon$ by the empirical distribution function of estimated residuals $\widetilde{\varepsilon}_t$, that is,

$$\widetilde{F}_\varepsilon(y) = \frac{1}{n} \sum_{t=1}^{n} I(\widetilde{\varepsilon}_t \leq y), \tag{3.1}$$

where $\widetilde{\varepsilon}_t = w(\mathbb{X}_t, X_t, \widehat{\theta})$. Note that the estimator $\widehat{\theta}$ used is assumed to satisfy (A2). Note, further, that instead of $\widetilde{F}_\varepsilon$, we can also use an estimator of $F_\varepsilon$ based on centered and/or standardized residuals like those used in the bootstrap schemes discussed in Section 4. This, however, will result in an extra term in the asymptotic covariances of the finite-dimensional distributions of the process $U_n$ given below.

The above considerations and the resulting specification of $F_{H_0}(\cdot|\mathbb{X}_{t-1})$ in (1.1) lead to the basic deviation process

$$U_n(x, y) = \frac{1}{\sqrt{n}} \sum_{t=1}^{n} I(\mathbb{X}_{t-1} \preceq x)[I(X_t \leq y) - \widetilde{F}_\varepsilon(w(\mathbb{X}_{t-1}, y, \widehat{\theta}))], \tag{3.2}$$

which is used in the sequel for testing the null hypothesis of interest. A basis for a test of $H_0$ is now given by the supremum deviation,

$$S_n = \sup_{(x,y) \in \mathbb{R}^p \times \mathbb{R}} |U_n(x, y)|. \tag{3.3}$$

Notice that, for any $x_1, \ldots, x_{k-1}, x_{k+1}, \ldots, x_p, y$ ($1 \leq k \leq p$), $U_n(x_1, \ldots, x_{k-1}, x_{k+1}, \ldots, x_p, y)$ is piecewise constant with possible jumps at $X_{1-k}, \ldots, X_{n-k}$. Furthermore, $U_n(x, \cdot)$ has possible jumps at $X_1, \ldots, X_n$ and is monotonously nonincreasing between these jumps. Hence, it follows that

$$S_n = \max_{(x,y) \in \mathcal{X}^n} |U_n(x, y)|,$$

where $\mathcal{X}^n = \{-\infty, X_0, \ldots, X_{n-1}\} \times \cdots \times \{-\infty, X_{1-p}, \ldots, X_{n-p}\} \times \{-\infty, X_1 - 0, X_1, \ldots, X_n - 0, X_n, \infty\}$; that is, it suffices to compute the test statistic by evaluating $U_n$ on the grid $\mathcal{X}^n$.

To derive the asymptotic distribution of $S_n$, we first determine the limit distribution of the processes $U_n$. Let $D = D(\bar{\mathbb{R}}^{p+1})$ be the space of cadlag functions on the extended $(p+1)$-dimensional Euclidean space $\bar{\mathbb{R}}^{p+1}$, that is, of functions which are continuous from above and possess limits from below. It is clear that $U_n$ belongs to $D$ with probability 1. Since we deal with suprema of these processes, it is convenient to endow $D$ with the supremum norm $\|\cdot\|$ and to prove weak convergence of the distributions in the normed space $(D, \|\cdot\|)$. In accordance with the discussion in Section IV.1 in Pollard (1984), we do not endow the space $D$ with the Borel $\sigma$-field generated by the closed sets under the uniform metric since this $\sigma$-field would be too rich, consequently creating measurability problems. Rather, we use the projection $\sigma$-field $\mathcal{P}$ generated by the coordinate projection maps. Since $S_n$ can be written as $\sup_{(x,y)\in\mathbb{Q}^p\times\mathbb{Q}}|U_n(x,y)|$, it is clear that it is $(\mathcal{P}-\mathcal{B})$-measurable, where $\mathcal{B}$ is the Borel $\sigma$-field. The fact that $D$ is not separable does not matter in the following since the limit process $U$ is concentrated on a separable subset of continuous functions. Consequently, we can apply the continuous mapping theorem (Theorem V.1 in Pollard (1984)) to derive the limit distribution of the test statistics. The following theorem establishes the asymptotic behavior of $S_n$ if the null hypothesis is true.

**Theorem 3.1.** *Suppose that $H_0$ is true with $X_t = G(\mathbb{X}_{t-1}, \theta, \varepsilon_t)$ in $\mathcal{M}_0$ satisfying* (A1) *or* (A1$'$). *Assume, further, that* (A2) *is fulfilled. If $n \to \infty$, then*

$$U_n \xrightarrow{d} U,$$

*where $U$ is a Gaussian process with continuous sample paths, zero mean and covariance function*

$$\begin{aligned}
&\Gamma((x_1, y_1), (x_2, y_2)) \\
&= \sum_{t\in\mathbb{Z}} \operatorname{cov}(g_1(\mathbb{X}_0, X_0; x_1, y_1) + g_2(\mathbb{X}_0, X_0; x_1, y_1) + g_3(\mathbb{X}_0, X_0; x_1, y_1), \\
&\qquad\qquad g_1(\mathbb{X}_{t-1}, X_t; x_2, y_2) + g_2(\mathbb{X}_{t-1}, X_t; x_2, y_2) + g_3(\mathbb{X}_{t-1}, X_t; x_2, y_2))
\end{aligned}$$

*with*

$$g_1(\mathbb{X}_{t-1}, X_t; x, y) = I(\mathbb{X}_{t-1} \preceq x)[I(X_t \leq y) - F_\varepsilon(w(\mathbb{X}_{t-1}, y, \theta))],$$

$$g_2(\mathbb{X}_{t-1}, X_t; x, y) = l(\mathbb{X}_{t-1}, X_t; \theta)' \int_{\{z \preceq x\}} \dot{F}_\varepsilon(w(z, y, \theta)) P^{\mathbb{X}_{t-1}}(\mathrm{d}z)$$

*and*

$$g_3(\mathbb{X}_{t-1}, X_t; x, y) = -\int_{\{z \preceq x\}} [I(\varepsilon_t \leq w(z, y, \theta)) - P(\varepsilon_t \leq w(z, y, \theta))] P^{\mathbb{X}_{t-1}}(\mathrm{d}z).$$

Theorem 3.1 and an application of the continuous mapping theorem together yield that an asymptotically $\alpha$-level test for testing $H_0$ is given by the following rule: reject $H_0$ if $S_n > t_{1-\alpha,\infty}$, where $t_{1-\alpha,\infty}$ denotes the $(1-\alpha)$-quantile point of the distribution of $\sup_{(x,y)\in\mathbb{R}^p\times\mathbb{R}}|U(x,y)|$. A bootstrap approach to calculating these quantiles is given in the next section.

We conclude this section by mentioning that testing more specific null hypotheses concerning the transition kernel of the underlying Markov process is also possible using our approach. For this, appropriate specifications of the basic deviation process (1.1) can be used, depending on the specifications imposed of $F_{H_0}(\cdot \,|\, \mathbb{X}_{t-1})$. Clearly, the most simple case is that of a fully specified Markov model, that is, that of testing

$$H_0^{(1)}: P_{\mathbf{X}} \in \{P_X : X_t | \mathbb{X}_{t-1} \sim F_0(\cdot \,|\, \mathbb{X}_{t-1})\} \cap \mathcal{M}.$$

The obvious specification of $U_n^{(0)}(x,y)$ which is appropriate in this case is

$$U_n^{(1)}(x,y) = \frac{1}{\sqrt{n}} \sum_{t=1}^{n} I(\mathbb{X}_{t-1} \preceq x)[I(X_t \leq y) - F_0(y \,|\, \mathbb{X}_{t-1})]. \tag{3.4}$$

Another example, which is of more interest in applications, is where the one-step transition distribution function under the null depends on an unknown vector of parameters $\theta$, but the distribution function of the innovations is specified, that is, where

$$H_0^{(2)}: P_{\mathbf{X}} \in \{P_X : X_t = G(\mathbb{X}_{t-1}, \theta, \varepsilon_t), \varepsilon_t \sim F_\varepsilon \text{ i.i.d.}; \theta \in \Theta\} \cap \mathcal{M}.$$

In this case and instead of (3.2), the specification

$$U_n^{(2)}(x,y) = \frac{1}{\sqrt{n}} \sum_{t=1}^{n} I(\mathbb{X}_{t-1} \preceq x)[I(X_t \leq y) - F_\varepsilon(w(\mathbb{X}_{t-1}, y, \widehat{\theta}))] \tag{3.5}$$

should be used. It can be shown, along the same lines as in the proof of Theorem 3.1, that $U_n^{(2)} \xrightarrow{d} U^{(2)}$, where $U^{(2)}$ is a zero-mean Gaussian process, the covariance function of which is that obtained from $\Gamma((x_1,y_1),(x_2,y_2))$ given in Theorem 3.1 after ignoring the component $g_3$.

## 4. Bootstrap approximations

### 4.1. The bootstrap procedure

To approximate the distribution of $S_n$ under the null hypotheses, we use a model-based bootstrap approach which employs the particular structure of the generating equation $X_t = G(\mathbb{X}_{t-1}, \theta, \varepsilon_t)$. In this context, the unknown parameter $\theta$ is replaced by its estimator $\widehat{\theta}$. Furthermore, the innovations $\varepsilon_t$ are replaced by pseudo-innovations generated according to the empirical distribution function of estimated errors. In particular, the pseudo-innovations are generated using the empirical distribution function

$$\widehat{F}_\varepsilon(y) = \frac{1}{n} \sum_{t=1}^{n} I(\widehat{\varepsilon}_t \leq y), \tag{4.1}$$

where, for instance, in the AR($p$)-case, $\widehat{\varepsilon}_t$ is given by $\widehat{\varepsilon}_t = \widetilde{\varepsilon}_t - n^{-1} \sum_{s=1}^{n} \widetilde{\varepsilon}_s$ and, in the ARCH($p$)-case, by $\widehat{\varepsilon}_t = \widetilde{\widetilde{\varepsilon}}_t / \sqrt{n^{-1} \sum_{s=1}^{n} \widetilde{\widetilde{\varepsilon}}_s^2}$ with $\widetilde{\widetilde{\varepsilon}}_t = \widetilde{\varepsilon}_t - n^{-1} \sum_{s=1}^{n} \widetilde{\varepsilon}_s$, where $\widetilde{\varepsilon}_s$ is defined in the sentence following equation (3.1).

The bootstrap algorithm used to approximate the distribution of $S_n$ under the null hypothesis is described by the following three steps.

1. Let $\mathbb{X}_0^* = (X_0^*, X_{-1}^*, \ldots, X_{1-p}^*)$ be some starting values. Given $\mathbb{X}_{t-1}^* = (X_{t-1}^*, X_{t-2}^*, \ldots, X_{t-p}^*)$, generate $X_t^*$ by

$$X_t^* = G(\mathbb{X}_{t-1}^*, \widehat{\theta}, \varepsilon_t^*),$$

   where $\varepsilon_t^*$ are i.i.d. random variables with $\varepsilon_t^* \sim \widehat{F}_\varepsilon$, where $\widehat{F}_\varepsilon$ is defined in (4.1).

2. Based on the bootstrap pseudo-series $(X_t^*)_{t=1-p,2-p,\ldots,n}$, let $U_n^*(x,y)$ be defined as $U_n(x,y)$ and obtained by replacing $\widehat{\theta}$ and $\widetilde{F}_\varepsilon$ in $U_n(x,y)$ by $\widehat{\theta}^*$ and $\widetilde{F}_\varepsilon^*$, respectively. Here, $\widehat{\theta}^*$ denotes the same estimator as $\widehat{\theta}$ based on $(X_t^*)_{t=1-p,2-p,\ldots,n}$. Furthermore,

$$\widetilde{F}_\varepsilon^*(y) = \frac{1}{n}\sum_{t=1}^n I(\widetilde{\varepsilon}_t^* \leq y), \tag{4.2}$$

   where $\widetilde{\varepsilon}_t^* = w(\mathbb{X}_t^*, X_t^*, \widehat{\theta}^*)$.

   The bootstrap analogue of $S_n$ is now given by

$$S_n^* = \sup_{(x,y) \in \mathbb{R}^p \times \mathbb{R}} |U_n^*(x,y)|.$$

3. Reject $H_0$ if

$$S_n > t_{1-\alpha,\infty}^*,$$

   where $t_{1-\alpha,\infty}^*$ denotes the $(1-\alpha)$-quantile of the distribution of $S_n^*$.

### 4.2. Some basic properties of the bootstrap processes

The derivation of theoretical results for the bootstrap is based on a case-by-case investigation since to establish properties such as stationarity and weak dependence of the bootstrap process, the particular model structure generating the $X_t^*$'s is explicitly used. Note that the bootstrap counterpart to $(X_t)_{t \in \mathbb{Z}}$ is the stationary (if it exists) process $(X_t^*)_{t \in \mathbb{Z}}$ obeying the equations

$$X_t^* = G(\mathbb{X}_{t-1}^*, \widehat{\theta}, \varepsilon_t^*). \tag{4.3}$$

We will see below that, in the AR($p$) and ARCH($p$) cases considered here, a unique solution to (4.3) exists with a probability tending to 1.

We first deal with the properties of the proposed bootstrap procedure to approximate the hypothesized conditional distributions $P^{X_t | \mathbb{X}_{t-1}}$ under the null hypothesis. The following lemma shows that the conditional distributions of the bootstrap process converge to the conditional distributions of the original process under the null and to some legitimate conditional distribution under the alternative. This convergence takes place in probability. To formulate such results in a transparent way, we define the following metric between distributions on $(\mathbb{R}^d, \mathcal{B}^d)$:

$$d(P,Q) = \inf_{X \sim P, Y \sim Q} \mathrm{E}[\|X-Y\| \wedge 1],$$

where the infimum is taken over all pairs $(X,Y)$ with $X \sim P$ and $Y \sim Q$, and where $\|\cdot\|$ is any norm on $\mathbb{R}^d$. A metric similar to this (the Mallows metric) has been used by Bickel and Freedman (1981), also in the context of proving bootstrap consistency. Convergence in the above metric is, in particular, equivalent to weak convergence. Concerning the behavior of the estimator $\widehat{\theta}$ in the case that the null hypothesis is not true, we make the following assumption.

(B1) There exists $\overline{\theta} \in \Theta$ such that $\widehat{\theta} \longrightarrow \overline{\theta}$ in probability.

**Lemma 4.1.** *Suppose that* (A1) *or* (A1') *and* (B1) *are fulfilled. Assume, further, that the density $f_e = F'_e$ of $e_t$ is bounded, where $e_t = w(\mathbb{X}_{t-1}, X_t, \overline{\theta})$. Then, for every compact set $K \subset \mathbb{R}^p$,*

$$\sup_{x \in K} d(P^{X_t^* | \mathbb{X}_{t-1}^* = x}, P_{\overline{\theta}}^{X_t | \mathbb{X}_{t-1} = x}) \xrightarrow{P} 0.$$

According to the above lemma, the asymptotic limit of $P^{X_t^* | \mathbb{X}_{t-1}^*}$ depends on whether the underlying hypothesis is true or not. In particular, $P^{X_t^* | \mathbb{X}_{t-1}^*}$ converges to that of the Markov process generated by $X_t = G(\mathbb{X}_{t-1}, \overline{\theta}, e_t)$, where the innovation sequence has distribution $F_e$. If $H_0$ is true, then $\overline{\theta} = \theta$ and $F_e = F_\varepsilon$. Note, however, that $P_{\mathbf{X}^*} \in \mathcal{M}_0$, where $P_{\mathbf{X}^*}$ denotes the law of $\mathbf{X}^* = \{X_t^*, t \in \mathbb{Z}\}$, even if the null hypothesis is not true. This is important for a good power behavior of the bootstrap-based test.

The following lemma is the key step in proving consistency properties of Markov bootstrap. It basically states that convergence of the conditional distributions implies convergence of the stationary distributions.

**Lemma 4.2.** *Suppose that $(Y_t)_{t \in \mathbb{Z}}$ and $(Y_t^{(n)})_{t \in \mathbb{Z}}$, $n \in \mathbb{N}$, are stationary Markov processes of order $p$, defined on probability spaces $(\Omega, \mathcal{A}, P)$ and $(\Omega_n, \mathcal{A}_n, P_n)$, respectively. Suppose that*

(i) *for all compact sets $K \subseteq \mathbb{R}^p$,*

$$\sup_{y \in K} d(P_n^{Y_t^{(n)} | Y_{t-1}^{(n)} = y_1, \ldots, Y_{t-p}^{(n)} = y_p}, P^{Y_t | Y_{t-1} = y_1, \ldots, Y_{t-p} = y_p}) \xrightarrow[n \to \infty]{} 0;$$

(ii) *for all $y \in \mathbb{R}^p$,*

$$\sup_{\widetilde{y}: \|\widetilde{y} - y\| \leq \delta} d(P^{Y_t | Y_{t-1} = y_1, \ldots, Y_{t-p} = y_p}, P^{Y_t | Y_{t-1} = \widetilde{y}_1, \ldots, Y_{t-p} = \widetilde{y}_p}) \xrightarrow[\delta \to 0]{} 0;$$

(iii) $(P_n^{Y_t^{(n)}})_{n \in \mathbb{N}}$ *is tight;*

(iv) *there is a unique stationary distribution $P^{Y_1,\ldots,Y_p}$ that corresponds to $P^{Y_t|Y_{t-1},\ldots,Y_{t-p}}$.*

*Then, for all $k \in \mathbb{N}$,*

$$P_n^{Y_1^{(n)},\ldots,Y_k^{(n)}} \Longrightarrow P^{Y_1,\ldots,Y_k}. \tag{4.4}$$

By Lemmas 4.1 and 4.2, we obtain the following result which shows the convergence of the finite-dimensional distributions of the bootstrap process to the desired joint distributions under the corresponding null hypothesis.

**Corollary 4.1.** *Suppose that the assumptions of Lemma 4.1 are satisfied. Then, for all $k \in \mathbb{N}$,*

$$d(P^{X_t^*,\ldots,X_{t+k}^*}, P_{\overline{\theta}}^{X_t,\ldots,X_{t+k}}) \xrightarrow{P} 0,$$

*where $P_{\overline{\theta}}^{X_t,\ldots,X_{t+k}}$ denotes the stationary probability measure of the Markov process $(X_t)_{t \in \mathbb{Z}}$ generated by $X_t = G(\mathbb{X}_{t-1}, \overline{\theta}, e_t)$ and where the i.i.d. innovation sequence satisfies $e_t \sim F_e$.*

Our next result deals with the weak dependence properties of the bootstrap processes.

**Lemma 4.3.** *Suppose that the assumptions of Lemma 4.1 are satisfied. There then exist sets $\Omega_n \subseteq \mathbb{R}^{n+p}$ such that $P((X_{1-p},\ldots,X_n) \in \Omega_n) \underset{n \to \infty}{\longrightarrow} 1$ and, for any sequence $(\omega_n)_{n \in \mathbb{N}}$ with $\omega_n \in \Omega_n$, $(X_t^*)_{t \in \mathbb{Z}}$ satisfies (conditionally under $(X_{1-p},\ldots,X_n)' = \omega_n$) conditions of weak dependence analogous to (2.2)–(2.4) and (2.6)–(2.8), respectively, with coefficients of weak dependence that can be majorized by a geometrically decaying series.*

### 4.3. Bootstrap validity

Based on the basic properties of the bootstrap procedure stated in the previous section, we are now able to justify asymptotically its use in obtaining critical values of the test statistics $S_n$. As in (A2), we assume the following.

(B2) The sequence of estimators $\widehat{\theta}^*$ admits the expansion

$$\widehat{\theta}^* - \widehat{\theta} = \frac{1}{n} \sum_{t=1}^{n} l(\mathbb{X}_{t-1}^*, X_t^*; \widehat{\theta}) + o_P(n^{-1/2}),$$

where $l(\cdot; \cdot)$ satisfies $\mathrm{E}^* l(\mathbb{X}_{t-1}^*, X_t^*; \widehat{\theta}) = 0$ and $\mathrm{E}^* \|l(\mathbb{X}_{t-1}^*, X_t^*; \widehat{\theta})\|^2 = O_P(1)$.

The next theorem establishes the asymptotic limit of $U_n^*$ which is used to evaluate the distribution of the test statistic $S_n$ under the null hypothesis.

**Theorem 4.1.** *Assume that* (A1) *or* (A1′) *as well as* (B1) *and* (B2) *are fulfilled. Suppose, further, that* $f_e = F'_e$ *is continuous and* $f_e(z) \underset{z \to \pm\infty}{\longrightarrow} 0$, *where* $F_e$ *denotes the distribution function of* $e_t = w(\mathbb{X}_{t-1}, X_t, \overline{\theta})$. *There then exist sets* $\widetilde{\Omega}_n \subset \mathbb{R}^{p+n}$ *such that* $P((X_{1-p}, \ldots, X_n)' \in \widetilde{\Omega}_n) \underset{n \to \infty}{\longrightarrow} 1$ *and, for every sequence* $(\omega_n)_{n \in \mathbb{N}}$ *with* $\omega_n \in \widetilde{\Omega}_n$, *we have (the bootstrap distribution is taken conditionally under* $(X_{1-p}, \ldots, X_n)' = \omega_n)$

$$U_n^* \xrightarrow{d} \overline{U}$$

*as* $n \to \infty$, *where* $\overline{U}$ *is a Gaussian process with continuous sample paths, zero mean and covariance function*

$$\overline{\Gamma}((x_1, y_1), (x_2, y_2))$$
$$= \sum_{t \in \mathbb{Z}} \mathrm{cov}(\overline{g}_1(\mathbb{X}_0, X_0; x_1, y_1) + \overline{g}_2(\mathbb{X}_0, X_0; x_1, y_1) + \overline{g}_3(\mathbb{X}_0, X_0; x_1, y_1),$$
$$\overline{g}_1(\mathbb{X}_{t-1}, X_t; x_2, y_2) + \overline{g}_2(\mathbb{X}_{t-1}, X_t; x_2, y_2) + \overline{g}_3(\mathbb{X}_{t-1}, X_t; x_2, y_2)),$$

*and* $\overline{g}_i(\mathbb{X}_{t-1}, X_t; x, y)$ *is defined for* $i = 1, 2, 3$ *as* $g_i(\mathbb{X}_{t-1}, X_t; x, y)$ *in Theorem* 3.1, *with* $\theta$ *replaced by* $\overline{\theta}$ *and* $F_\varepsilon$ *by* $F_e$.

Note that if $H_0$ is true, then $\widetilde{U} = U$ since $\overline{\theta} = \theta$ and $F_e = F_\varepsilon$. In this case, the limiting behavior of the bootstrap statistic $S_n^*$ is identical to that of the statistics $S_n$ given in Theorem 3.1. On the other hand, if $H_0$ is not true, then $t_{1-\alpha, \infty}^* \to c$ as $n \to \infty$, where $c$ denotes the $(1-\alpha)$-quantile point of the limiting distribution of $\sup_{(x,y) \in \mathbb{R}^p \times \mathbb{R}} |\overline{U}(x, y)|$. In this case, $n^{-1/2} S_n \to C$ in probability, where $C$ denotes a positive constant. Therefore,

$$\lim_{n \to \infty} P(S_n > t_{1-\alpha, \infty}^*) = \begin{cases} \alpha, & \text{if } H_0 \text{ is true,} \\ 1, & \text{if } H_1 \text{ is true,} \end{cases}$$

that is, the test based on the bootstrap critical values $t_{1-\alpha, \infty}^{*(i)}$ asymptotically achieves the desired level $\alpha$ and is consistent.

We conjecture that our test has nontrivial power for local alternatives converging to the null at a $\sqrt{n}$-rate. To illustrate this, consider the simple case of a fully specified null hypothesis $H_0^{(1)}$ and sequences of local alternatives corresponding to Markov processes having one-step transition distribution functions given by

$$F_n(\cdot | \mathbb{X}_{t-1}) = F_0(\cdot | \mathbb{X}_{t-1}) + \frac{1}{\sqrt{n}} H(\cdot, \mathbb{X}_{t-1}),$$

where $H(\cdot)$ is an appropriate function satisfying

$$\frac{1}{n} \sum_{t=1}^{n} I(\mathbb{X}_{t-1} \preceq x) H(y, \mathbb{X}_{t-1}) \to D(y, x) \not\equiv 0 \quad (4.5)$$

in probability as $n \to \infty$. In this case, the corresponding deviation process $U_n^{(1)}(y,x)$ can be decomposed as

$$U_n^{(1)}(y,x) = \frac{1}{\sqrt{n}} \sum_{t=1}^n I(\mathbb{X}_{t-1} \preceq x)[I(X_t \leq y) - F_n(y|\mathbb{X}_{t-1})]$$
$$+ \frac{1}{n} \sum_{t=1}^n I(\mathbb{X}_{t-1} \preceq x) H(y, \mathbb{X}_{t-1}),$$

from which (taking into account (4.5)) the desired result will follow by showing that

$$\widetilde{U}_n^{(1)}(y,x) = \frac{1}{\sqrt{n}} \sum_{t=1}^n I(\mathbb{X}_{t-1} \preceq x)[I(X_t \leq y) - F_n(y|\mathbb{X}_{t-1})]$$

converges to the same Gaussian process as the process (3.4).

## 5. Numerical examples

**Example 1.** The test statistic $S_n^{(2)} = \max_{(x,y) \in \mathcal{X}^n} |U_n^{(2)}(x,y)|$ is applied to test the hypothesis that the underlying process obeys the ARCH(1) structure

$$X_t = \sqrt{\theta_0 + \theta_1 X_{t-1}^2} \varepsilon_t$$

with independent and standard Gaussian-distributed errors. Three different sample sizes, $n = 100$, 200 and 400, have been considered. The parameters of the process have been set equal to $\theta_0 = 0.1$ and $\theta_1 = 0.4$. The results obtained are based on the least-squares estimator of $\theta_0$ and $\theta_1$. Table 1 presents the empirical rejection probabilities. To investigate the power of our test procedure, different types of alternatives have been considered. One alternative to the hypothesis of an ARCH(1) process with Gaussian errors is where the distribution of the innovations is given by $\varepsilon_t = \eta_t \sqrt{(\nu - 2)/\nu}$, where $(\eta_t)_{t \in \mathbb{Z}}$ is an i.i.d. sequence of $t$-distributed random variables with $\nu$ degrees of freedom. We also investigated the power of our test for the case where the true process is an ARCH(2) process with additional parameter $\theta_2 = 0.4$, a GARCH(1, 1) process $X_t = \sigma_t \varepsilon_t$ with $\sigma_t^2 = 0.08 + 0.7 X_{t-1}^2 + 0.2 \sigma_{t-1}^2$ and the simple stochastic volatility (SV) model $X_t = \exp\{h_t/2\} \varepsilon_t$ and $h_t = -0.9 + 0.6 h_{t-1} + \omega_t$, where $\omega_t$ is a sequence of independent standard Gaussian random variables.

**Table 1.** Empirical rejection probabilities for testing the hypothesis of an ARCH(1) model

|  | $\alpha$ | $H_0^{(2)}$ TRUE: ARCH(1) $\varepsilon_t \sim \mathcal{N}(0,1)$ | $H_0^{(2)}$ FALSE: ARCH(1) $\varepsilon_t \sim t_5$ | ARCH(2) $\theta_2 = 0.4$ | GARCH(1, 1) | SV-Model |
|---|---|---|---|---|---|---|
| $n = 100$ | 0.05 | 0.064 | 0.202 | 0.190 | 0.151 | 0.460 |
|  | 0.10 | 0.109 | 0.310 | 0.278 | 0.235 | 0.660 |
| $n = 200$ | 0.05 | 0.041 | 0.281 | 0.295 | 0.225 | 0.775 |
|  | 0.10 | 0.115 | 0.415 | 0.402 | 0.345 | 0.871 |
| $n = 400$ | 0.05 | 0.046 | 0.485 | 0.515 | 0.366 | 0.980 |
|  | 0.10 | 0.112 | 0.671 | 0.635 | 0.495 | 0.991 |

***Example 2.*** Using the test statistic $S_n$, the hypothesis of interest is that the underlying process is an i.i.d. process with standard Gaussian- or uniform-distributed innovations. The alternative considered to this null hypothesis is that the underlying process is a first order autoregressive process $X_t = \theta X_{t-1} + \varepsilon_t$ with $\varepsilon_t \sim \mathcal{N}(0, 1)$ and three different values of $\theta$. Table 2 presents the empirical rejection probabilities of $S_n$ for testing the corresponding null hypothesis based on sample sizes of length $n = 25$, $n = 50$ and $n = 100$.

The results presented in Tables 1 and 2 are based on 200 replications of the underlying process, where, for each replication, critical values of the test have been obtained using 500 bootstrap samples. Although computational requirements prevented us from considering larger sample sizes, more trials, more bootstrap replications or more complicated models, the results obtained are very encouraging. In particular, and as these tables show, the test statistic retains the desired size under the null hypothesis and shows a very good power behavior for the different types of alternatives considered.

***A real-data example.*** We apply our testing procedure to the first $n = 2000$ observations of the monthly log-returns of the Intel stock series analyzed in Tsay (2005). Tsay (2005), page 109, selected, for this series, the ARCH(1) model

$$r_t = 0.0174 + X_t, \qquad X_t = \sqrt{0.0134 + 0.2492 X_{t-1}^2} \varepsilon_t,$$

with standard Gaussian-distributed innovations $\varepsilon_t$. For this model, the value of the test statistic $S_n^{(2)} = \max_{(x,y) \in \mathcal{X}^n} |U_n^{(2)}(x, y)|$ equals 15.025, while a bootstrap estimate of the upper 5% percentage point of the distribution of the same statistic under the null equals

**Table 2.** Empirical rejection probabilities for testing the hypothesis of an i.i.d. sequence

|  |  | $H_0^{(3)}$ TRUE: $\varepsilon_t \sim \mathcal{N}(0,1)$ | $H_0^{(3)}$ FALSE: $\varepsilon_t \sim U(-\sqrt{3}, \sqrt{3})$ | $\theta = 0.2$ | $\theta = 0.4$ | $\theta = 0.6$ |
| --- | --- | --- | --- | --- | --- | --- |
| $n = 25$ | 0.05 | 0.075 | 0.081 | 0.105 | 0.330 | 0.625 |
|  | 0.10 | 0.152 | 0.161 | 0.268 | 0.485 | 0.740 |
| $n = 50$ | 0.05 | 0.041 | 0.058 | 0.170 | 0.585 | 0.889 |
|  | 0.10 | 0.104 | 0.105 | 0.270 | 0.709 | 0.925 |
| $n = 100$ | 0.05 | 0.062 | 0.069 | 0.305 | 0.834 | 0.995 |
|  | 0.10 | 0.115 | 0.118 | 0.455 | 0.925 | 0.998 |

1.501. This percentage point has been estimated using $B = 1000$ bootstrap replications. Our testing procedure therefore leads to a rejection of the above ARCH(1) model for the Intel stock series. Note that standard methods, based on residuals, for checking the fit of the above model do not indicate any inadequacy of the fitted ARCH(1) model in describing the conditional heteroscedasticity of the data; see Tsay (2005), page 111, for details.

## 6. A central limit theorem

The first central limit theorems for weakly dependent sequences were given by Corollary A in Doukhan and Louhichi (1999) and Theorem 1 in Coulon-Prieur and Doukhan (2000). While the former result is for sequences of stationary random variables, the latter one is tailor-made for triangular arrays of asymptotically sparse random variables as they appear with kernel density estimators. Below, we state a central limit theorem for general triangular schemes of weakly dependent random variables. An interesting aspect of this result is that no moment condition beyond Lindeberg's is required.

**Theorem 6.1.** *Suppose that* $(X_{n,k})_{k=1,\ldots,n}$, $n \in \mathbb{N}$, *is a triangular scheme of (row-wise) stationary random variables with* $\mathrm{E} X_{n,k} = 0$ *and* $\mathrm{E} X_{n,k}^2 \leq C < \infty$. *Furthermore, we assume that*

$$\frac{1}{n} \sum_{k=1}^{n} \mathrm{E} X_{n,k}^2 I(|X_{n,k}|/\sqrt{n} > \epsilon) \underset{n \to \infty}{\longrightarrow} 0 \tag{6.1}$$

*holds for all* $\epsilon > 0$ *and that*

$$\mathrm{var}(X_{n,1} + \cdots + X_{n,n})/n \underset{n \to \infty}{\longrightarrow} \sigma^2 \in [0, \infty). \tag{6.2}$$

*For* $n \geq n_0$, *there exists a monotonously nonincreasing and summable sequence* $(\theta_r)_{r \in \mathbb{N}}$ *such that, for all indices* $1 \leq s_1 < s_2 < \cdots < s_u < s_u + r = t_1 \leq t_2 \leq n$, *the following upper bounds for covariances hold true: for all measurable and quadratic integrable functions* $f : \mathbb{R}^u \longrightarrow \mathbb{R}$,

$$|\mathrm{cov}(f(X_{n,s_1}, \ldots, X_{n,s_u}), X_{n,t_1})| \leq \sqrt{\mathrm{E} f^2(X_{n,s_1}, \ldots, X_{n,s_u})} \theta_r \tag{6.3}$$

*and for all measurable and bounded functions* $f : \mathbb{R}^u \longrightarrow \mathbb{R}$,

$$|\mathrm{cov}(f(X_{n,s_1}, \ldots, X_{n,s_u}), X_{n,t_1} X_{n,t_2})| \leq \|f\|_\infty \theta_r, \tag{6.4}$$

*where* $\|f\|_\infty = \sup_{x \in \mathbb{R}^u} |f(x)|$. *Then,*

$$\frac{1}{\sqrt{n}}(X_{n,1} + \cdots + X_{n,n}) \overset{d}{\longrightarrow} \mathcal{N}(0, \sigma^2).$$

**Proof.** If $\sigma^2 = 0$, then we obviously have that $\frac{1}{\sqrt{n}}(X_{n,1}+\cdots+X_{n,n}) \xrightarrow{d} \mathcal{N}(0,0)$. Therefore, it remains to prove the assertion in the case $\sigma^2 > 0$, which we assume in the rest of the proof. Let $Y_{n,k} = X_{n,k}/\sqrt{\mathrm{E}(X_{n,1}+\cdots+X_{n,n})^2}$. In view of (6.2), it suffices to show that

$$Y_{n,1} + \cdots + Y_{n,n} \xrightarrow{d} \mathcal{N}(0,1). \tag{6.5}$$

To prove this, we use the classical Lindeberg method which was first adapted to causal CLT's by Rio (1995).

We set $\sigma_n^2 = \mathrm{var}(X_{n,1}+\cdots+X_{n,n}) = n\,\mathrm{var}(X_{n,1}) + 2\sum_{j=1}^{n-1}(n-j)\mathrm{cov}(X_{n,1},X_{n,j+1})$ and $v_{n,k} = \mathrm{var}(Y_{n,1}+\cdots+Y_{n,k}) - \mathrm{var}(Y_{n,1}+\cdots+Y_{n,k-1}) = (\mathrm{var}(X_{n,1}) + 2\sum_{j=1}^{k-1}\mathrm{cov}(X_{n,1},X_{n,j+1}))/\sigma_n^2$. We obtain, by $|\mathrm{cov}(X_{n,1},X_{n,j+1})| \le C\theta_j$, that

$$|nv_{n,k} - 1| \le \frac{2C}{\sigma_n^2}\left(\sum_{j=1}^{k-1} j\theta_j + \sum_{j=k}^{n-1}(n-j)\theta_j\right).$$

Since, by majorized convergence, $\sum_{j=1}^{\infty}(j/n)\theta_j \xrightarrow[n\to\infty]{} 0$ and $\sum_{j=k}^{\infty}\theta_j \xrightarrow[k\to\infty]{} 0$, it follows that there exist $k_0, n_0 \in \mathbb{N}$ such that

$$v_{n,k} \ge 0 \qquad \text{for all } (n,k) \text{ with } n \ge n_0 \text{ and } k_0 \le k \le n. \tag{6.6}$$

To simplify the notation in the rest of the proof, we pretend that (6.6) holds for all $(n,k)$ with $n \ge 1$ and $1 \le k \le n$. (Otherwise, we start with $n_0$ and sum the first $k_0 - 1$ random variables in each row to a new random variable, $Y_{n,0} = Y_{n,1}+\cdots+Y_{n,k_0-1}$. We then prove the assertion for the sums $Y_{n,0}+Y_{n,k_0}+\cdots+Y_{n,n}$.)

Let $h:\mathbb{R} \longrightarrow \mathbb{R}$ be an arbitrary, three times continuously differentiable function with $\|h^{(j)}\|_\infty =: C_j < \infty$, $j = 0, \ldots, 3$. Furthermore, let $Z_{n,k} \sim \mathcal{N}(0, v_{n,k})$, $k = 1, \ldots, n$, be independent random variables which are also independent of $(Y_{n,k})_{k=1,\ldots,n}$. Since $v_{n,1}+\cdots+v_{n,n} = 1$, it follows from Theorem 7.1 in Billingsley (1968) that it suffices to show that

$$\mathrm{E}h(Y_{n,1}+\cdots+Y_{n,n}) - \mathrm{E}h(Z_{n,1}+\cdots+Z_{n,n}) \xrightarrow[n\to\infty]{} 0. \tag{6.7}$$

We define $S_{n,k} = \sum_{j=1}^{k-1} Y_{n,j}$ and $T_{n,k} = \sum_{j=k+1}^{n} Z_{n,j}$. Then,

$$\mathrm{E}h(Y_{n,1}+\cdots+Y_{n,n}) - \mathrm{E}h(Z_{n,1}+\cdots+Z_{n,n}) = \sum_{k=1}^{n} \Delta_{n,k},$$

where

$$\Delta_{n,k} = \mathrm{E}[h(S_{n,k}+Y_{n,k}+T_{n,k}) - h(S_{n,k}+Z_{n,k}+T_{n,k})].$$

We further decompose $\Delta_{n,k} = \Delta_{n,k}^{(1)} - \Delta_{n,k}^{(2)}$, where

$$\Delta_{n,k}^{(1)} = \mathrm{E}h(S_{n,k}+Y_{n,k}+T_{n,k}) - \mathrm{E}h(S_{n,k}+T_{n,k}) - \frac{v_{n,k}}{2}\mathrm{E}h^{(2)}(S_{n,k}+T_{n,k}),$$

$$\Delta_{n,k}^{(2)} = \mathrm{E}h(S_{n,k} + Z_{n,k} + T_{n,k}) - \mathrm{E}h(S_{n,k} + T_{n,k}) - \frac{v_{n,k}}{2}\mathrm{E}h^{(2)}(S_{n,k} + T_{n,k}).$$

We will show that

$$\sum_{k=1}^{n} \Delta_{n,k}^{(i)} \xrightarrow[n\to\infty]{} 0 \qquad \text{for } i = 1, 2.$$

(i) Upper bound for $|\sum_{k=1}^{n} \Delta_{n,k}^{(2)}|$.
Since $\mathrm{E}Z_{n,k}h'(S_{n,k} + T_{n,k}) = 0$, we have, for some random $\rho_{n,k} \in (0,1)$, that

$$\Delta_{n,k}^{(2)} = \mathrm{E}\frac{Z_{n,k}^2}{2}[h^{(2)}(S_{n,k} + \rho_{n,k}Z_{n,k} + T_{n,k}) - h^{(2)}(S_{n,k} + T_{n,k})].$$

Hence, we obtain that

$$\left|\sum_{k=1}^{n} \Delta_{n,k}^{(2)}\right| \leq \frac{C_3}{2} \sum_{k=1}^{n} \mathrm{E}|Z_{n,k}|^3$$
$$\leq \frac{C_3}{2}\mathrm{E}|\mathcal{N}(0,1)|^3 \max_{1\leq k\leq n}\{\sqrt{v_{n,k}}\} \xrightarrow[n\to\infty]{} 0. \tag{6.8}$$

(ii) Upper bound for $|\sum_{k=1}^{n} \Delta_{n,k}^{(1)}|$.
Let $\epsilon > 0$ be arbitrary. We will actually show that

$$\left|\sum_{k=1}^{n} \Delta_{n,k}^{(1)}\right| \leq \epsilon \qquad \text{for all } n \geq n(\epsilon). \tag{6.9}$$

We have, for some random $\tau_{n,k} \in (0,1)$, that

$$\Delta_{n,k}^{(1)} = \mathrm{E}Y_{n,k}h'(S_{n,k} + T_{n,k}) + \mathrm{E}\left[\frac{Y_{n,k}^2}{2}h^{(2)}(S_{n,k} + \tau_{n,k}Y_{n,k} + T_{n,k})\right]$$
$$- \frac{v_{n,k}}{2}\mathrm{E}h^{(2)}(S_{n,k} + T_{n,k}).$$

Since $\mathrm{E}Y_{n,k}h'(T_{n,k}) = 0$, we have, again for some random $\mu_{n,k,j} \in (0,1)$, that

$$\mathrm{E}Y_{n,k}h'(S_{n,k} + T_{n,k}) = \sum_{j=1}^{k-1} \mathrm{E}Y_{n,k}[h'(S_{n,j+1} + T_{n,k}) - h'(S_{n,j} + T_{n,k})]$$
$$= \sum_{j=1}^{k-1} \mathrm{E}Y_{n,k}Y_{n,j}h^{(2)}(S_{n,j} + \mu_{n,k,j}Y_{n,j} + T_{n,k}).$$

This yields, in conjunction with

$$v_{n,k} = EY_{n,k}^2 + 2\sum_{j=1}^{k-1} EY_{n,k}Y_{n,j},$$

that

$$\Delta_{n,k}^{(1)} = \sum_{j=1}^{k-d} EY_{n,k}Y_{n,j}[h^{(2)}(S_{n,j} + \mu_{n,k,j}Y_{n,j} + T_{n,k}) - Eh^{(2)}(S_{n,k} + T_{n,k})]$$

$$+ \sum_{j=k-d+1}^{k-1} EY_{n,k}Y_{n,j}[h^{(2)}(S_{n,j} + \mu_{n,k,j}Y_{n,j} + T_{n,k}) - Eh^{(2)}(S_{n,k} + T_{n,k})]$$

$$+ \tfrac{1}{2}EY_{n,k}^2[h^{(2)}(S_{n,k} + \tau_{n,k}Y_{n,k} + T_{n,k}) - Eh^{(2)}(S_{n,k} + T_{n,k})]$$

$$= \Delta_{n,k}^{(1,1)} + \Delta_{n,k}^{(1,2)} + \Delta_{n,k}^{(1,3)},$$

say. (The value of $d$ does *not* depend on $n$ and its proper choice is indicated below.)

We now have, by (6.3), that

$$|\Delta_{n,k}^{(1,1)}| \leq \sum_{j=1}^{k-d} \sqrt{EY_{n,k}^2} O\left(\frac{1}{\sqrt{n}}\right)\theta_{k-j} = O\left(n^{-1}\sum_{j=d}^{n}\theta_j\right).$$

By choosing $d$ sufficiently large, we obtain that

$$\left|\sum_{k=1}^{n}\Delta_{n,k}^{(1,1)}\right| \leq \frac{\epsilon}{3} \qquad \text{for all } n \geq n(\epsilon). \tag{6.10}$$

The term $\Delta_{n,k}^{(1,2)}$ will be split up as

$$\Delta_{n,k}^{(1,2)} = \sum_{j=k-d+1}^{k-1} EY_{n,k}Y_{n,j}[h^{(2)}(S_{n,j} + \mu_{n,k,j}Y_{n,j} + T_{n,k}) - h^{(2)}(S_{n,j-d} + T_{n,k})]$$

$$+ \sum_{j=k-d+1}^{k-1} EY_{n,k}Y_{n,j}[h^{(2)}(S_{n,j-d} + T_{n,k}) - Eh^{(2)}(S_{n,j-d} + T_{n,k})]$$

$$+ \sum_{j=k-d+1}^{k-1} EY_{n,k}Y_{n,j}[Eh^{(2)}(S_{n,j-d} + T_{n,k}) - Eh^{(2)}(S_{n,k} + T_{n,k})]$$

$$= \Delta_{n,k}^{(1,2,1)} + \Delta_{n,k}^{(1,2,2)} + \Delta_{n,k}^{(1,2,3)},$$

say. The Lindeberg condition (6.1) yields that, for arbitrary $\epsilon' > 0$,

$$\left| \sum_{k=1}^{n} \Delta_{n,k}^{(1,2,1)} \right|$$
$$\leq 2C_2 \sqrt{\sum_{k=1}^{n} \mathrm{E}Y_{n,k}^2 I(|Y_{n,k}| > \epsilon')} \sqrt{\sum_{k=1}^{n} \sum_{j=k-d+1}^{k-1} \mathrm{E}Y_{n,j}^2} + \epsilon' \sqrt{\sum_{k=1}^{n} \sum_{j=k-d+1}^{k-1} \mathrm{E}Y_{n,j}^2}$$
$$\times \sqrt{\sum_{k=1}^{n} \sum_{j=k-d+1}^{k-1} \mathrm{E}[h^{(2)}(S_{n,j} + \mu_{n,k,j} Y_{n,j} + T_{n,k}) - h^{(2)}(S_{n,j-d} + T_{n,k})]^2}$$
$$= o(1) + O(\epsilon').$$

Using condition (6.4), we obtain that

$$|\Delta_{n,k}^{(1,2,2)}| = O(n^{-1} \mathrm{d}\theta_d).$$

From the monotonicity and summability of the sequence $(\theta_k)_{k \in \mathbb{N}}$, it follows that $\mathrm{d}\theta_d \xrightarrow[d \to \infty]{} 0$. Furthermore, the relation

$$\Delta_{n,k}^{(1,2,3)} = O(n^{-3/2})$$

is obvious. Again, for sufficiently large $d$, these upper estimates yield that

$$\left| \sum_{k=1}^{n} \Delta_{n,k}^{(1,2)} \right| \leq \frac{\epsilon}{3} \qquad \text{for all } n \geq n(\epsilon). \tag{6.11}$$

Finally, we obtain, in complete analogy to the calculations above, that

$$\sum_{k=1}^{n} \Delta_{n,k}^{(1,3)} \xrightarrow[n \to \infty]{} 0, \tag{6.12}$$

which completes, in conjunction with (6.10) and (6.11), the proof of (6.9). □

## 7. Proofs of auxiliary lemmas and main results

Proofs of some of our main results are given in this section, while, for some others, we stress only the essentials. More details, as well as the proofs of Lemma 2.1 and Theorem 4.1, which are omitted in the sequel, are given in Neumann and Paparoditis (2005).

There has been much effort made in the literature to prove stochastic equicontinuity, often as a sufficient condition for tightness, of families of multivariate processes. For a family of processes $(X_n)_{n\in\mathbb{N}}$ with sample paths in $C([0,1]^q)$, one seeks to show that, for all $\delta > 0$, $\eta > 0$, there exists an $\varepsilon > 0$ such that

$$P(\omega_\varepsilon(X_n) > \delta) \leq \eta \qquad \forall n \geq n_0, \tag{7.1}$$

where the modulus of continuity is defined as $\omega_\varepsilon(x) = \sup_{\|s-t\|\leq\varepsilon} |x(s) - x(t)|$. For a family of processes on $\mathbb{R}^q$, one can either transform them to processes on $[0,1]^q$ or one can alternatively show that, for all $\delta > 0$, $\eta > 0$, there exists a grid $\mathcal{G} = \{b_0^{(1)}, \ldots, b_M^{(1)}\} \times \cdots \times \{b_0^{(q)}, \ldots, b_M^{(q)}\}$ with $-\infty = b_0^{(r)} < \cdots < b_M^{(r)} = \infty$, $r = 1, \ldots, q$, such that

$$P\left(\max_{1\leq i_1,\ldots,i_q\leq M} \sup_{b(\underline{i}-\mathbb{1}_q)\preceq t\preceq b(\underline{i})} |X_n(t) - X_n(b(\underline{i}))| > \delta\right) \leq \eta, \tag{7.2}$$

for all $n \geq n_0$. Here, and in the following, we use the notation $\underline{i} = (i_1, \ldots, i_q)$, $t = (t_1, \ldots, t_q)$ and $b(\underline{i}) = (b_{i_1}^{(1)}, \ldots, b_{i_q}^{(q)})$.

It has been shown, for example, in Pollard (1984), Theorem 3 in Section V.1, that (7.2) and the weak convergence of finite-dimensional distributions of $X_n$ to those of a process $X$ with continuous sample paths together imply that $(X_n)_{n\in\mathbb{N}}$ converges in distribution (with respect to the supremum norm) to $X$. (Pollard actually proves this for processes on $[0,1]$; the extension to processes on $\mathbb{R}^q$ is, however, obvious.) An obstacle to proving (7.2) arises since the supremum over an infinite set is involved. Therefore, one often proves, instead of (7.2), that there exists a sequence of increasingly fine grids, $\mathcal{G}_n = \{t_{n,1}^{(1)}, \ldots, t_{n,M_n}^{(1)}\} \times \cdots \times \{t_{n,1}^{(q)}, \ldots, t_{n,M_n}^{(q)}\}$, such that $\mathcal{G} \subseteq \mathcal{G}_n$ and

$$P\left(\max_{1\leq i_1,\ldots,i_q\leq M} \sup_{t\in\mathcal{G}_n: b(\underline{i}-\mathbb{1}_q)\preceq t\preceq b(\underline{i})} |X_n(t) - X_n(b(\underline{i}))| > \delta\right) \leq \eta \tag{7.3}$$

holds for $n \geq n_0$, and then derives (7.2) by continuity and monotonicity arguments. With a slight abuse of terminology, we also call property (7.2) stochastic equicontinuity and (7.3) stochastic equicontinuity on the grids $\mathcal{G}_n$.

In the following, we provide a simple proof of (7.3) based on an anisotropic dyadic tiling of the space. Such an anisotropic tiling has previously been used in the proof of Proposition 7.3 in Rio (2000), page 100ff. This proof constitutes an alternative to the commonly used approach based on Bickel and Wichura's (1971) fluctuation result for their modulus of continuity $M''$, and to an approach based on an isotropic tiling of the space proposed by Neuhaus (1971).

**Lemma 7.1.** *Let $(X_n(t))_{t\in\mathbb{R}^q}$ be a sequence of real-valued stochastic processes. For any hyperrectangle $B = (s_1, t_1] \times \cdots \times (s_q, t_q]$, the increment of $X_n$ around $B$ is given by*

$$X_n(B) = \sum_{(\varepsilon_1,\ldots,\varepsilon_q)\in\{0,1\}^q} (-1)^{q-(\varepsilon_1+\cdots+\varepsilon_q)} X_n(s_1 + \varepsilon_1(t_1-s_1), \ldots, s_q + \varepsilon_q(t_q-s_q)).$$

*We suppose that there exists a sequence of measures $(\mu_n)_{n\in\mathbb{N}}$ on $(\mathbb{R}^q, \mathcal{B}^q)$ with continuous marginals which converges weakly to a finite measure $\mu$, also having continuous marginals,*

*and*

$$\mathrm{E}[X_n(B)]^4 \leq [\mu_n(B) + Cn^{-q}]^{1+\gamma} \tag{7.4}$$

*for all hyperrectangles $B$ and some $\gamma > 0$, $C < \infty$. Let $\delta$ and $\eta$ be arbitrary positive constants.*

*There then exists a coarse grid $\mathcal{G} = \{b_0^{(1)}, \ldots, b_M^{(1)}\} \times \cdots \times \{b_0^{(q)}, \ldots, b_M^{(q)}\}$ with $-\infty = b_0^{(r)} < \cdots < b_M^{(r)} = \infty$ ($r = 1, \ldots, q$) and a sequence of fine grids $\mathcal{G}_n = \{t_{n,1}^{(1)}, \ldots, t_{n,M_n}^{(1)}\} \times \cdots \times \{t_{n,1}^{(q)}, \ldots, t_{n,M_n}^{(q)}\}$ with $\mathcal{G} \subseteq \mathcal{G}_n$ such that*

$$\mu_n(\mathbb{R}^{r-1} \times (t_{n,k-1}^{(r)}, t_{n,k}^{(r)}] \times \mathbb{R}^{q-r}) \leq 2\mu(\mathbb{R}^q) n^{-1} \qquad \forall k = 1, \ldots, M_n, \forall r = 1, \ldots, q \tag{7.5}$$

*and*

$$P\left(\max_{1 \leq i_1, \ldots, i_q \leq M} \max_{t \in \mathcal{G}_n : b(\underline{i} - \mathbb{1}_q) \preceq t \preceq b(\underline{i})} |X_n(t) - X_n(b(\underline{i}))| > \delta\right) \leq \eta \tag{7.6}$$

*holds for all $n \geq n_0$, where $n_0$ is sufficiently large.*

**Remark 7.1.** As already mentioned in the discussion after Theorem 3 of Bickel and Wichura (1971), it is possible (and, for the bootstrap processes, important) that the measures $\mu_n$ are allowed to depend on $n$. The term $n^{-q}$ in (7.4) cannot be avoided in our context; see, for example, (7.15) in the proof of Theorem 3.1 below. Because of this term, we obtain stochastic equicontinuity in a first step only on a grid with cardinality of $M_n^q = O(n^q)$. In our applications, stochastic equicontinuity over the full space will then follow from monotonicity and continuity properties of the processes involved; see step (i) in the proof of Theorem 3.1 below.

**Proof of Lemma 7.1.** (i) *Dyadic systems of grid points.* First, we define dyadic systems of grid points. At the coarse scales, their choice is tied to the measure $\mu$. Let $F^{(r)}$ be the $r$th marginal cumulative distribution function of $\mu$, that is,

$$F^{(r)}(x) = \mu(\mathbb{R}^{r-1} \times (-\infty, x] \times \mathbb{R}^{q-r}).$$

For an appropriate $J_0 \in \mathbb{N}$ to be determined in part (iii) of this proof, we define, for $r = 1, \ldots, q$ and $0 \leq j \leq J_0$,

$$b_{j,k}^{(r)} = \begin{cases} -\infty, & \text{if } k = 0, \\ F^{(r)^{-1}}(k 2^{-j} \mu(\mathbb{R}^q)), & \text{if } 1 \leq k < 2^j, \\ \infty, & \text{if } k = 2^j. \end{cases}$$

We choose $M = 2^{J_0}$ and $b_k^{(r)} = b_{J_0, k}^{(r)}$.

At the finer scales, with index $j > J_0$, the grid points are chosen according to $\mu_n$. We set

$$F_n^{(r)}(x) = \mu_n(\mathbb{R}^{r-1} \times (-\infty, x] \times \mathbb{R}^{q-r})$$

and choose, for $j = J_0 + 1, \ldots, J_n$ with $2^{J_n-1} < n \le 2^{J_n}$, grid points as follows. For $l \in \{0, \ldots, 2^{J_0} - 1\}$ and $k \in \{1, \ldots, 2^{j-J_0}\}$, we define $b^{(r)}_{j,l2^{j-J_0}+k} = b^{(r)}_{j,l2^{j-J_0}+k}(n)$ such that

$$F_n^{(r)}(b^{(r)}_{j,l2^{j-J_0}+k}(n)) = F_n^{(r)}(b_{J_0,l}) + \frac{k}{2^{j-J_0}}[F_n^{(r)}(b_{J_0,l+1}) - F_n^{(r)}(b_{J_0,l})].$$

(Again, we set $b^{(r)}_{j,0}(n) = -\infty$ and $b^{(r)}_{j,2^j}(n) = \infty$.) We set $M_n = 2^{J_n} + 1$ and $t^{(r)}_{n,1} = b^{(r)}_{J_n,0}, \ldots, t^{(r)}_{n,M_n} = b^{(r)}_{J_n,2^{J_n}}$ ($r = 1, \ldots, q$). That is, the fine grids are given as $\mathcal{G}_n = \mathcal{G}_n^{(1)} \times \cdots \times \mathcal{G}_n^{(q)}$, where $\mathcal{G}_n^{(r)} = \{b^{(r)}_{J_n,0}, \ldots, b^{(r)}_{J_n,2^{J_n}}\}$.

Since $\mu_n \Longrightarrow \mu$, we have, for all $r = 1, \ldots, q$, $l = 1, \ldots, 2^{J_0}$ and $n \ge n_0$ with $n_0$ sufficiently large, that

$$\mu_n(\mathbb{R}^{r-1} \times (b^{(r)}_{J_0,l-1}, b^{(r)}_{J_0,l}] \times \mathbb{R}^{q-r}) \le 2\mu(\mathbb{R}^{r-1} \times (b^{(r)}_{J_0,l-1}, b^{(r)}_{J_0,l}] \times \mathbb{R}^{q-r}) = 2^{1-J_0}\mu(\mathbb{R}^q).$$

This implies that

$$\mu_n(\mathbb{R}^{r-1} \times (b^{(r)}_{j,k-1}, b^{(r)}_{j,k}] \times \mathbb{R}^{q-r}) \le 2^{1-j}\mu(\mathbb{R}^q) \tag{7.7}$$

for all $r = 1, \ldots, q$, $j = J_0 + 1, \ldots, J_n$, $k = 1, \ldots, 2^j$ and $n \ge n_0$, that is, (7.5) is satisfied.

(ii) *A probabilistic bound for the increments of $X_n$.* To simplify notation, in the sequel, we use multiindices $\underline{j} = (j_1, \ldots, j_q)$ and $\underline{k} = (k_1, \ldots, k_q)$. For $(\underline{j}, \underline{k})$ from the set $\mathcal{B}_n = \{(\underline{j}, \underline{k}) : 0 \le j_r \le 2^{J_n}, 1 \le k_r \le 2^{j_r} \ \forall r\}$, we define the hyperrectangle

$$B_{\underline{j},\underline{k}} = (b^{(1)}_{j_1,k_1-1}, b^{(1)}_{j_1,k_1}] \times \cdots \times (b^{(q)}_{j_q,k_q-1}, b^{(q)}_{j_q,k_q}].$$

We choose any $\alpha \in (0, \gamma/4)$ and define the thresholds

$$\lambda_{\underline{j}} = K 2^{-\alpha(j_1 + \cdots + j_q)/q}, \tag{7.8}$$

where $K$ will be chosen below. From (7.7) and $2^{J_n} \le 2n$ we obtain that

$$\mu_n(B_{\underline{j},\underline{k}}) + Cn^{-q} \le C_1 \min_{1 \le r \le q}\{2^{-j_r}\} \le C_1 2^{-(j_1 + \cdots + j_q)/q} \qquad \forall (\underline{j}, \underline{k}) \in \mathcal{B}_n.$$

Therefore, we obtain, by Markov's inequality, that

$$P(|X_n(B_{\underline{j},\underline{k}})| > \lambda_{\underline{j}}) \le \frac{[\mu_n(B_{\underline{j},\underline{k}}) + n^{-q}]^{1+\gamma}}{K^4 \cdot 2^{-4\alpha(j_1 + \cdots + j_q)/q}}$$

$$\le (\mu_n(B_{\underline{j},\underline{k}}) + n^{-q})K^{-4}C_1^{\gamma}2^{(j_1 + \cdots + j_q)(4\alpha - \gamma)/q}.$$

This implies that

$$P(|X_n(B_{\underline{j},\underline{k}})| > \lambda_{\underline{j}} \text{ for any } (\underline{j}, \underline{k}) \in \mathcal{B}_n)$$

$$\leq (\mu_n(\mathbb{R}^q) + 2^{J_n q} n^{-q}) K^{-4} C_1^\gamma \sum_{j_1,\ldots,j_q=0}^{J_n} 2^{(j_1+\cdots+j_q)(4\alpha-\gamma)/q} \qquad (7.9)$$

$$\leq \eta,$$

for $n \geq n_0$, provided the constant $K$ in (7.8) is large enough.

(iii) *Stochastic equicontinuity of $X_n$ on the fine grid.* Now, assume that

$$|X_n(B_{\underline{j},\underline{k}})| \leq \lambda_{\underline{j}} \qquad \text{for all } (\underline{j},\underline{k}) \in \mathcal{B}_n. \qquad (7.10)$$

Moreover, let $t \in \mathcal{G}_n$ and $\underline{i}$ be such that $b(\underline{i}-\mathbb{1}_q) \preceq t \preceq b(\underline{i})$. There then exist hyperrectangles $B_{\underline{j}^{(1)},\underline{k}^{(1)}}, \ldots, B_{\underline{j}^{(L)},\underline{k}^{(L)}}$ with *different scale indices* $\underline{j}^{(l)}$ and $(\underline{j}^{(l)},\underline{k}^{(l)}) \in \mathcal{B}_n$ such that $\max_{1\leq r\leq q} j_r^{(l)} \geq J_0$ and

$$\{x : x \preceq b(\underline{i})\} \setminus \{x : x \preceq t\} = \bigcup_{l=1}^{L} B_{\underline{j}^{(l)},\underline{k}^{(l)}}. \qquad (7.11)$$

Accordingly, $X_n(b(\underline{i})) - X_n(t) = \sum_{l=1}^{L} X_n(B_{\underline{j}^{(l)},\underline{k}^{(l)}})$, which implies that

$$|X_n(t) - X_n(b(\underline{i}))| \leq \sum_{l=1}^{L} \lambda_{\underline{j}^{(l)}} \leq Kq \sum_{j_1=J_0}^{J_n} \sum_{j_2,\ldots,j_q=0}^{J_n} 2^{-\alpha(j_1+\cdots+j_q)/q}. \qquad (7.12)$$

Now, choosing $J_0 \in \mathbb{N}$ such that $Kq2^{-J_0\alpha/q}(1/(1-2^{-\alpha/q}))^q \leq \delta$, we obtain that

$$|X_n(t) - X_n(b(\underline{i}))| \leq \delta$$

holds for all $\underline{i}$ and $t \in \mathcal{G}_n$ with $b(\underline{i}-\mathbb{1}_q) \preceq t \preceq b(\underline{i})$, whenever (7.10) is fulfilled. This, however, implies, in conjunction with (7.9), that (7.6) is satisfied. □

**Proof of Theorem 3.1.** We apply the method of proving weak convergence for processes described by, for example, Wichura (1971), Proposition 1 and Pollard (1984), Theorem 3 in Section V.1. To this end, we will prove (i) stochastic equicontinuity of $(U_n)_{n\in\mathbb{N}}$ and (ii) weak convergence of the finite-dimensional distributions. From step (i), we can identify the prospective limit process $U$ as a centered Gaussian one. Because of the complicated covariance function, we cannot immediately see that $U$ possesses a version with continuous sample paths. However, steps (i) and (ii) together imply that there is a version of $U$ which inherits the property of stochastic continuity from the processes $(U_n)_{n\in\mathbb{N}}$. Having this, it is then easy to conclude that this process has continuous sample paths with probability 1. These facts together yield the desired convergence of $(U_n)_{n\in\mathbb{N}}$ to $U$.

(i) *Stochastic equicontinuity of $(U_n)_{n\in\mathbb{N}}$.* We set $q = p+1$. We will show that there exists, for any $\delta > 0$ and $\eta > 0$, a grid $\mathcal{G} = \{b_0^{(1)},\ldots,b_M^{(1)}\} \times \cdots \times \{b_0^{(q)},\ldots,b_M^{(q)}\}$ with

Goodness-of-fit tests for Markovian time series models 39

$$-\infty = b_0^{(r)} < b_1^{(r)} < \cdots < b_M^{(r)} = \infty \ (r=1,\ldots,q) \text{ such that (with } b(\underline{i}) = (b_{i_1}^{(1)},\ldots,b_{i_q}^{(q)})')$$

$$P\left(\max_{1 \le i_1,\ldots,i_q \le M} \sup_{t \in \mathbb{R}^q : b(\underline{i} - \mathbb{1}_q) \preceq t \preceq b(\underline{i})} |U_n(t) - U_n(b(\underline{i}))| > \delta\right) \le \eta. \quad (7.13)$$

We decompose $U_n(x,y)$ as

$$\begin{aligned}
U_n(x,y) &= \frac{1}{\sqrt{n}} \sum_{t=1}^n I(\mathbb{X}_{t-1} \preceq x)[I(X_t \le y) - F_\varepsilon(w(\mathbb{X}_{t-1}, y, \theta))] \\
&\quad + \frac{1}{\sqrt{n}} \sum_{t=1}^n I(\mathbb{X}_{t-1} \preceq x)[F_\varepsilon(w(\mathbb{X}_{t-1}, y, \theta)) - F_\varepsilon(w(\mathbb{X}_{t-1}, y, \widehat{\theta}))] \\
&\quad + \frac{1}{\sqrt{n}} \sum_{t=1}^n I(\mathbb{X}_{t-1} \preceq x)[F_\varepsilon(w(\mathbb{X}_{t-1}, y, \widehat{\theta})) - \widehat{F}_\varepsilon(w(\mathbb{X}_{t-1}, y, \widehat{\theta}))] \\
&=: R_n^{(1)}(x,y) + R_n^{(2)}(x,y) + R_n^{(3)}(x,y),
\end{aligned} \quad (7.14)$$

say. It is now most convenient to prove stochastic equicontinuity for $R_n^{(1)}$, $R_n^{(2)}$ and $R_n^{(3)}$ separately. We will give all details for $R_n^{(1)}$ and refer to Neumann and Paparoditis (2005) for more details regarding $R_n^{(2)}$ and $R_n^{(3)}$. For any hyperrectangle $B = B_x \times B_y = (s_1, t_1] \times \cdots \times (s_q, t_q]$, denote by

$$R_n^{(1)}(B) = \frac{1}{\sqrt{n}} \sum_{t=1}^n I(\mathbb{X}_{t-1} \in B_x)[I(X_t \in B_y) - P(X_t \in B_y \mid \mathbb{X}_{t-1})]$$

the increment of $R_n^{(1)}$ around $B$. We will first show that, for all $\gamma \in (0, 1/q)$, there exists some constant $C_\gamma < \infty$ such that

$$\mathrm{E}[R_n^{(1)}(B)]^4 \le C_\gamma [P^{\mathbb{X}_{t-1}, X_t}(B) + n^{-q}]^{1+\gamma}. \quad (7.15)$$

Let $g_B(\mathbb{X}_{t-1}, X_t) = I(\mathbb{X}_{t-1} \in B_x)[I(X_t \in B_y) - P(X_t \in B_y \mid \mathbb{X}_{t-1})]$. Note that, for sufficiently integrable random variables $Y_1, \ldots, Y_4$ with $EY_i = 0$, the relations

$$\begin{aligned}
EY_1 \cdots Y_4 &= EY_1Y_2 \cdot EY_3Y_4 + \mathrm{cov}(Y_1Y_2, Y_3Y_4) \\
&= \mathrm{cov}(Y_1, Y_2Y_3Y_4) \\
&= \mathrm{cov}(Y_1Y_2Y_3, Y_4)
\end{aligned}$$

hold true. Since $g_B$ is a bounded function with $\mathrm{E}g_B(\mathbb{X}_{t-1}, X_t) = 0$, we consequently obtain that

$$\mathrm{E}[R_n^{(1)}(B)]^4 \le \frac{4!}{n^2} \sum_{t_1 \le t_2 \le t_3 \le t_4} |\mathrm{E}[g_B(\mathbb{X}_{t_1-1}, X_{t_1}) \cdots g_B(\mathbb{X}_{t_4-1}, X_{t_4})]|$$

$$\leq 4! \left\{ \left[ \sum_{r=0}^{n-1} C_{r,2}(B) \right]^2 + \frac{3}{n} \sum_{r=0}^{n-1} (r+1)^2 C_{r,4}(B) \right\}, \tag{7.16}$$

where

$$C_{r,q}(B) = \max_{1 \leq m \leq q-1} \sup_{(t_1,\ldots,t_q) \in \mathcal{T}_{r,q}(m)} |\operatorname{cov}(g_B(\mathbb{X}_{t_1-1}, X_{t_1}) \cdots g_B(\mathbb{X}_{t_m-1}, X_{t_m}),$$
$$g_B(\mathbb{X}_{t_{m+1}-1}, X_{t_{m+1}}) \cdots g_B(\mathbb{X}_{t_q-1}, X_{t_q}))|$$

and

$$\mathcal{T}_{r,q}(m) = \left\{ (t_1,\ldots,t_q) : 1 \leq t_1 \leq \ldots \leq t_q \leq n, \max_{1 \leq j \leq q-1} \{t_{j+1} - t_j\} = t_m - t_{m+1} = r \right\}.$$

(Inequality (7.16) is similar to inequality (2.14) in Doukhan and Louhichi (1999), the only difference being that the '3' is absent there.)

For small values of $r$, we use the simple estimate

$$C_{r,q}(B) \leq \operatorname{E}|g_B(\mathbb{X}_{t-1}, X_t)| \leq 2P^{\mathbb{X}_{t-1}, X_t}(B). \tag{7.17}$$

For large values of $r$, we intend to exploit the weak dependence of the process $(X_t)_{t \in \mathbb{Z}}$ in order to show that $C_{r,q}(B)$ gets small as $r$ increases. Since $g_B$ is not Lipschitz, we define smooth approximations to $g_B$,

$$g_{B,\epsilon}(x,y) = \int w_\epsilon(u) g_{B+u}(x,y) \, \mathrm{d}u,$$

where $(w_\epsilon)_{\epsilon > 0}$ is a family of nonnegative functions with $\operatorname{supp}(w_\epsilon) \subseteq \{u = (u_1,\ldots,u_{p+1}) : u_i \geq 0 \text{ and } \|u\|_{l_1} \leq \epsilon\}$, $\int w_\epsilon(u) \, \mathrm{d}u = 1$ and $\|w_\epsilon\|_\infty \leq C\epsilon^{-q}$. Since $\operatorname{Lip}(g_{B,\epsilon}(u_1) \cdots g_{B,\epsilon}(u_m)) \leq m \cdot \operatorname{Lip}(g_{B,\epsilon}) = O(1/\epsilon)$, we obtain, by (2.2) or (2.6), respectively, that

$$|\operatorname{cov}(g_{B,\epsilon}(\mathbb{X}_{t_1-1}, X_{t_1}) \cdots g_{B,\epsilon}(\mathbb{X}_{t_m-1}, X_{t_m}),$$
$$g_{B,\epsilon}(\mathbb{X}_{t_{m+1}-1}, X_{t_{m+1}}) \cdots g_{B,\epsilon}(\mathbb{X}_{t_1-1}, X_{t_q}))| \tag{7.18}$$
$$\leq C \frac{\rho^{r-p}}{\epsilon}.$$

Since $|x_1 \cdots x_q - y_1 \cdots y_q| \leq \sum_{i=1}^q |x_i - y_i|$ for all real numbers $x_i, y_i \in [-1,1]$ and, by Lipschitz continuity of $F_X$,

$$\operatorname{E}|g_B(\mathbb{X}_{t-1}, X_t) - g_{B,\epsilon}(\mathbb{X}_{t-1}, X_t)| \leq C\epsilon,$$

we obtain that

$$|\operatorname{cov}(g_B(\mathbb{X}_{t_1-1}, X_{t_1}) \cdots g_B(\mathbb{X}_{t_m-1}, X_{t_m}),$$

$$g_B(\mathbb{X}_{t_{m+1}-1}, X_{t_{m+1}}) \cdots g_B(\mathbb{X}_{t_1-1}, X_{t_q}))$$
$$- \text{cov}(g_{B,\epsilon}(\mathbb{X}_{t_1-1}, X_{t_1}) \cdots g_{B,\epsilon}(\mathbb{X}_{t_m-1}, X_{t_m}), \qquad (7.19)$$
$$g_{B,\epsilon}(\mathbb{X}_{t_{m+1}-1}, X_{t_{m+1}}) \cdots g_{B,\epsilon}(\mathbb{X}_{t_q-1}, X_{t_q}))|$$
$$\leq C\epsilon.$$

From (7.17) and (7.18) and (7.19) with $\epsilon = \rho^{(r-p)/2}$, we obtain that

$$C_{r,q}(B) \leq C(P^{\mathbb{X}_{t-1}, X_t}(B) \wedge \rho^{(r-p)/2}),$$

which implies, by (7.16), inequality (7.15).

From (7.15), we conclude, by Lemma 7.1, that there exists a coarse grid $\widetilde{\mathcal{G}} = \{b_0^{(1)}, \ldots, b_M^{(1)}\} \times \cdots \times \{b_0^{(q)}, \ldots, b_M^{(q)}\}$ and a sequence of fine grids $\widetilde{\mathcal{G}}_n = \{t_{n,1}^{(1)}, \ldots, t_{n,M_n}^{(1)}\} \times \cdots \times \{t_{n,1}^{(q)}, \ldots, t_{n,M_n}^{(q)}\}$ with $\widetilde{\mathcal{G}} \subseteq \widetilde{\mathcal{G}}_n$ such that

$$F_X(t_{n,k}^{(r)}) - F_X(t_{n,k-1}^{(r)}) \leq C n^{-1}$$

and

$$P\left(\max_{1 \leq i_1, \ldots, i_q \leq M} \max_{t \in \widetilde{\mathcal{G}}_n : b(\underline{i}-\mathbb{1}_q) \preceq t \preceq b(\underline{i})} |R_n^{(1)}(t) - R_n^{(1)}(b(\underline{i}))| > \frac{\delta}{6}\right) \leq \frac{\eta}{6} \qquad (7.20)$$

for all $n \geq n_0$ and $n_0$ sufficiently large.

To extend property (7.20) to the whole space, we employ a simple monotonicity argument. For $t$ with $t_{n,i_r-1}^{(r)} \leq t_r \leq t_{n,i_r}^{(r)}$ $\forall r$, we have the inequalities

$$R_n^{(1)}(t_{n,i_1-1}^{(1)}, \ldots, t_{n,i_q-1}^{(q)})$$
$$- \frac{1}{\sqrt{n}} \sum_{t=1}^n [I(\mathbb{X}_{t-1} \preceq (t_{n,i_1}^{(1)}, \ldots, t_{n,i_p}^{(p)})) P(X_t \leq t_{n,i_{p+1}}^{(p+1)} \mid \mathbb{X}_{t-1})$$
$$- I(\mathbb{X}_{t-1} \preceq (t_{n,i_1-1}^{(1)}, \ldots, t_{n,i_p-1}^{(p)})) P(X_t \leq t_{n,i_{p+1}-1}^{(p+1)} \mid \mathbb{X}_{t-1})]$$
$$\leq R_n^{(1)}(t_1, \ldots, t_q)$$
$$\leq R_n^{(1)}(t_{n,i_1}^{(1)}, \ldots, t_{n,i_q}^{(q)}) + \frac{1}{\sqrt{n}} \sum_{t=1}^n [\cdots].$$

It follows from the Bernstein-type inequality for weakly dependent random variables, from Kallabis and Neumann (2006), that

$$P\left(\frac{1}{\sqrt{n}} \sum_{t=1}^n [\cdots] > \frac{\delta}{6} \text{ for any } (i_1, \ldots, i_q)\right) \leq \frac{\eta}{6}. \qquad (7.21)$$

(7.20) and (7.21) together yield that

$$P\left(\max_{1\leq i_1,\ldots,i_q\leq M}\max_{t\in\mathbb{R}^q:b(\underline{i}-\mathbb{1}_{p+1})\preceq t\preceq b(\underline{i})}|R_n^{(1)}(t)-R_n^{(1)}(b(\underline{i}))|>\frac{\delta}{3}\right)\leq\frac{\eta}{3} \quad (7.22)$$

for all $n\geq n_0$.

Furthermore, we can also prove property (7.22) for the processes $R_n^{(2)}$ and $R_n^{(3)}$, possibly with other coarse grids $\widehat{\mathcal{G}}$ and $\overline{\mathcal{G}}$. This yields property (7.13) for the grid $\mathcal{G}$, which is the combination of the grids $\widetilde{\mathcal{G}}$, $\widehat{\mathcal{G}}$ and $\overline{\mathcal{G}}$.

(ii) *Weak convergence of finite-dimensional distributions.* Let $k\in\mathbb{N}$, $(x_1,y_1),\ldots,(x_k,y_k)\in\mathbb{R}^p\times\mathbb{R}$ and $c_1,\ldots,c_k\in\mathbb{R}$ be arbitrary. By the Cramér–Wold device, it suffices to show that

$$\sum_{l=1}^{k}c_lU_n(x_l,y_l)\stackrel{d}{\longrightarrow}\mathcal{N}\left(0,\sum_{l,m=1}^{k}c_lc_m\Gamma((x_l,y_l),(x_m,y_m))\right). \quad (7.23)$$

According to (7.14), we can show that

$$U_n(x,y)=\frac{1}{\sqrt{n}}\sum_{t=1}^{n}[g_1(\mathbb{X}_{t-1},X_t;x,y)+\cdots+g_3(\mathbb{X}_{t-1},X_t;x,y)]+o_P(1).$$

Moreover, it follows from (2.2)–(2.4) and (2.6)–(2.8), respectively, that the triangular scheme $(Z_{n,t})_{t=1,\ldots,n}$ with $Z_{n,t}=\sum_{l=1}^{k}c_l[g_1(\mathbb{X}_{t-1},X_t;x_l,y_l)+\cdots+g_3(\mathbb{X}_{t-1},X_t;x_l,y_l)]$ satisfies conditions (6.1)–(6.4) of the central limit theorem in Section 6. (Note that the function $g_1(\cdot,\cdot;x,y)$ is discontinuous; we must use an approximation by a smoothed version, as above, for checking (6.2) to (6.4).) Hence, (7.23) follows immediately from Theorem 6.1.

(iii) *The limit process.* According to (7.13), there exists a sequence of grids $\mathcal{G}^{(N)}=\{b_{1,0}^{(N)},\ldots,b_{1,M_N}^{(N)}\}\times\cdots\times\{b_{q,0}^{(N)},\ldots,b_{q,M_N}^{(N)}\}$ with $-\infty=b_{r,0}^{(N)}<\cdots<b_{r,M_N}^{(N)}=\infty$ ($r=1,\ldots,q$) and $\mathcal{G}^{(N+1)}\subseteq\mathcal{G}^{(N)}$ $\forall N\in\mathbb{N}$ such that

$$P\left(\max_{1\leq i_1,\ldots,i_q\leq M_N}\sup_{t\in\mathbb{R}^q:b^{(N)}(\underline{i}-\mathbb{1}_q)\preceq t\preceq b^{(N)}(\underline{i})}|U_n(t)-U_n(b^{(N)}(\underline{i}))|\geq\frac{1}{N}\right)\leq\frac{1}{N} \quad (7.24)$$

for $n\geq n_N$ and $n_N$ sufficiently large.

Let $\mathcal{G}^{(\infty)}=\bigcup_{N=1}^{\infty}\mathcal{G}^{(N)}$. Kolmogorov's consistency theorem ensures that there exists a real-valued stochastic process $\widetilde{U}$ on $\mathcal{G}^{(\infty)}$ with finite-dimensional distributions

$$\begin{pmatrix}\widetilde{U}(t_1)\\ \vdots\\ \widetilde{U}(t_k)\end{pmatrix}\sim\mathcal{N}\left(0_k,\begin{pmatrix}\Gamma(t_1,t_1)&\cdots&\Gamma(t_1,t_k)\\ \vdots&\ddots&\vdots\\ \Gamma(t_k,t_1)&\cdots&\Gamma(t_k,t_k)\end{pmatrix}\right).$$

It follows from the weak convergence proved in step (ii) that $\widetilde{U}$ inherits the continuity property from the processes $(U_n)_{n\in\mathbb{N}}$, that is, we have that

$$P\left(\max_{1\leq i_1,\ldots,i_q\leq M_N}\max_{t\in\mathcal{G}(N'):b^{(N)}(\underline{i}-\mathbb{1}_q)\preceq t\preceq b^{(N)}(\underline{i})}|\widetilde{U}(t)-\widetilde{U}(b^{(N)}(\underline{i}))|\geq \frac{1}{N}\right) \tag{7.25}$$
$$\leq \frac{1}{N} \qquad \forall N'\geq N.$$

By continuity from below of the probability measure $P$, this implies that even

$$P\left(\max_{1\leq i_1,\ldots,i_q\leq M_N}\max_{t\in\mathcal{G}(\infty):b^{(N)}(\underline{i}-\mathbb{1}_q)\preceq t\preceq b^{(N)}(\underline{i})}|\widetilde{U}(t)-\widetilde{U}(b^{(N)}(\underline{i}))|\geq \frac{1}{N}\right)\leq \frac{1}{N} \tag{7.26}$$

holds true.

We now extend the process $\widetilde{U}$ to a process $U$ on $\mathbb{R}^q$. Fix any $t\in\mathbb{R}^q$. There then exists a sequence $(\underline{i}^{(N)})_{N\in\mathbb{N}}$ such that $b^{(N)}(\underline{i}^{(N)}-\mathbb{1}_q)\preceq t\preceq b^{(N)}(\underline{i}^{(N)})$. It follows from (7.26) that $\widetilde{U}(b^{(N)}(\underline{i}^{(N)}))$ converges in probability. We set

$$U(t)=\plim_{N\to\infty}\widetilde{U}(b^{(N)}(\underline{i}^{(N)})).$$

It is clear that the process $U$ has the same stochastic continuity property as $\widetilde{U}$, that is,

$$P\left(\max_{1\leq i_1,\ldots,i_q\leq M_N}\sup_{t\in\mathbb{R}^q:b^{(N)}(\underline{i}-\mathbb{1}_q)\preceq t\preceq b^{(N)}(\underline{i})}|U(t)-U(b^{(N)}(\underline{i}))|\geq \frac{1}{N}\right)\leq \frac{1}{N}. \tag{7.27}$$

This means, in particular, that $U$ has, with probability 1, continuous sample paths. Moreover, it follows from (7.13), (7.23) and (7.27) that the finite-dimensional distributions of $U_n$ converge to those of $U$. The latter property yields, in conjunction with (7.13) and (7.27), by Theorem 3 in Section V.1 of Pollard (1984), that $U_n\xrightarrow{d}U$. □

**Proof of Lemma 4.1.** First, note that for the function $w(\cdot)$, as it is specified under (A1) or (A1'), there exists $\delta>0$ such that, for all $\theta_1,\theta_2\in U_\delta(\overline{\theta})=\{\vartheta\in\Theta,\|\vartheta-\overline{\theta}\|<\delta\}$ and all compact sets $K\subset\mathbb{R}^{p+1}$, $|w(x,y,\theta_1)-w(x,y,\theta_2)|\leq L(x,y)\|\theta_1-\theta_2\|$ with $\sup_{x,y\in K}L(x,y)<\infty$. For such a $\delta$ and $\widetilde{K}<\infty$, define

$$\Omega_n=\{(X_{1-p},\ldots,X_n):\|\widehat{\theta}(X_{1-p},\ldots,X_n)-\overline{\theta}\|<\delta,\mathrm{E}((\varepsilon_t^*)^2\mid X_{1-p},\ldots,X_n)\leq \widetilde{K}\}.$$

Choose $\widetilde{K}<\infty$ large enough so that $P((X_{1-p},\ldots,X_n)\in\Omega_n)\to 1$ as $n\to\infty$. We then have, for $\omega_n\in\Omega_n$,

$$|P(X_t^*\leq y|\mathbb{X}_{t-1}^*=x)-F_e(w(x,y,\overline{\theta}))|$$
$$=|\widehat{F}_\varepsilon(w(x,y,\widehat{\theta}))-F_e(w(x,y,\overline{\theta}))|$$
$$\leq|F_e(w(x,y,\widehat{\theta}))-F_e(w(x,y,\overline{\theta}))|+|\widehat{F}_\varepsilon(w(x,y,\widehat{\theta}))-F_e(w(x,y,\widehat{\theta}))|.$$

For the first term on the right-hand side of the last inequality above, we obtain, using (B1), that

$$\sup_{x \in K} |F_e(w(x,y,\widehat{\theta})) - F_e(w(x,y,\overline{\theta}))| \le \|f_e\|_\infty \|\widehat{\theta} - \overline{\theta}\| \sup_{x \in K} L(x,y) \to 0$$

as $n \to \infty$. The second term can be majorized by $\|\widehat{F}_\varepsilon - F_e\|_\infty$, which converges to zero in probability; see Neumann and Paparoditis (2005) for details. □

**Proof of Lemma 4.2.** Let $P^{(k)} = P^{Y_1,\ldots,Y_k}$ and $P_n^{(k)} = P^{Y_1^{(n)},\ldots,Y_k^{(n)}}$. For $k = p$, the proof that (4.4) holds true follows the lines of the proof of Theorem 3.3 in Paparoditis and Politis (2002).

For $k < p$, the assertion follows from (4.4) with $k = p$, by the continuous mapping theorem. For $k > p$, the assertion follows by induction from the result for $k = p$ and the fact that, for any bounded and uniformly continuous function $f : \mathbb{R}^k \longrightarrow \mathbb{R}$, the relation

$$\begin{aligned}
\int f \, \mathrm{d}P_n^{(k)} &= \int_{\mathbb{R}^{k-1}} [\mathrm{E}(f(Y_1^{(n)},\ldots,Y_k^{(n)}) \mid (Y_1^{(n)},\ldots,Y_{k-1}^{(n)})' = y) \\
&\qquad - \mathrm{E}(f(Y_1,\ldots,Y_k) \mid (Y_1,\ldots,Y_{k-1})' = y)] P_n^{(k-1)}(\mathrm{d}y) \\
&\quad + \int_{\mathbb{R}^{k-1}} \mathrm{E}(f(Y_1,\ldots,Y_k) \mid (Y_1,\ldots,Y_{k-1})' = y) [P_n^{(k-1)}(\mathrm{d}y) - P^{(k-1)}(\mathrm{d}y)] \\
&\quad + \int f \, \mathrm{d}P^{(k)} \\
&\underset{n \to \infty}{\longrightarrow} \int f \, \mathrm{d}P^{(k)}
\end{aligned}$$

holds true. □

**Proof of Corollary 4.1.** For some null sequence $(\delta_n)_{n \in \mathbb{N}}$ and appropriate $K < \infty$, we define a set of "favorable events" such that

$$\widetilde{\Omega}_n \subseteq \{(X_{1-p},\ldots,X_n) : \|\widehat{\theta}_n(X_1,\ldots,X_n) - \overline{\theta}\| \le \delta_n, \\ \mathrm{E}((\varepsilon_t^*)^2 \mid X_{1-p},\ldots,X_n) \le K\}. \tag{7.28}$$

Moreover, let the $\widetilde{\Omega}_n$ be such that, for any sequence $(\omega_n)_{n \in \mathbb{N}}$ with $\omega_n \in \widetilde{\Omega}_n$,

$$\mathcal{L}(\varepsilon_t^* \mid (X_{1-p},\ldots,X_n) = \omega_n) \Longrightarrow \mathcal{L}(e_t).$$

The constant $K < \infty$ and the sequence $(\delta_n)_{n \in \mathbb{N}}$ above are chosen such that $\delta_n \to 0$ and $P((X_{1-p},\ldots,X_n) \in \widetilde{\Omega}_n) \to 1$ as $n \to \infty$. Now, let $(\omega_n)_{n \in \mathbb{N}}$ be an arbitrary sequence with $\omega_n \in \widetilde{\Omega}_n$. We now assume that the bootstrap distributions are taken under the condition that $(X_{1-p},\ldots,X_n) = \omega_n$. (This refers, in general, to a triangular scheme, but not to a single sequence of $X_t$.) Since $\widehat{\theta}_n(X_{1-p},\ldots,X_n) \underset{n\to\infty}{\longrightarrow} \overline{\theta}$ and $E((\varepsilon_t^*)^2 \mid X_{1-p},\ldots,X_n) \leq K$, the conditions of Lemma 4.2 are fulfilled, which yields the assertion. $\square$

**Proof of Lemma 4.3.** We only stress the essentials of the proof. Let $(\omega_n)_{n\in\mathbb{N}}$ be any sequence with $\omega_n \in \widetilde{\Omega}_n$, where $\widetilde{\Omega}_n$ is chosen as in the proof of Corollary 4.1. Assume that the distributions are taken under the condition $(X_{1-p},\ldots,X_n) = \omega_n$.

For the AR($p$) case, let $\xi_1,\ldots,\xi_p$ be the roots of the polynomial $\overline{\theta}(\cdot)$ with $\overline{\theta}(z) = 1 - \overline{\theta}_1 z - \cdots - \overline{\theta}_p z^p$. According to (A1), we have that $\epsilon := \min\{|\xi_1|,\ldots,|\xi_p|\} - 1 > 0$. If $\delta_n$ in (7.28) is sufficiently small, then we obtain, for the roots $\widehat{\xi}_{n,1},\ldots,\widehat{\xi}_{n,p}$ of $\widehat{\theta}(z) = 1 - \widehat{\theta}_{n,1} z - \cdots - \widehat{\theta}_{n,p} z^p$, that $\min\{|\widehat{\xi}_{n,1}|,\ldots,|\widehat{\xi}_{n,p}|\} \geq 1 + \epsilon/2$; see Theorem 1.4 in Marden (1949). Thus, there exists a stationary solution to the equation $X_t^* = \widehat{\theta}_{n,1} X_{t-1}^* + \cdots + \widehat{\theta}_{n,p} X_{t-p}^* + \varepsilon_t^*$ which can be written as an MA($\infty$)-process, $X_t^* = \sum_{k=0}^{\infty} \beta_{n,k} \varepsilon_{t-k}^*$, where $|\beta_{n,k}| \leq C_\delta (1 + \epsilon/2 - \delta)^{-k}$ for any $\delta > 0$ and corresponding $C_\delta < \infty$. The rest of the proof then follows that in Section 2.1.

To show the weak dependence of the bootstrap process in the ARCH($p$), choose $\delta_n < \eta/p$ in (7.28), where $\eta = 1 - \sum_{i=1}^p \overline{\theta}_i > 0$, and note that, for all $\widehat{\theta} = (\widehat{\theta}_1,\ldots,\widehat{\theta}_p) \in U_\delta(\overline{\theta})$, there exists some $\widetilde{\eta} \in (0,1)$ such that $\sum_{i=1}^p \widehat{\theta}_i < 1 - \widetilde{\eta}$. Furthermore, a (unique) stationary solution to the equation $X_t^* = \varepsilon_t^* \sqrt{\widehat{\theta}_0 + \widehat{\theta}_1 X_{t-1}^{*2} + \cdots + \widehat{\theta}_p X_{t-p}^{*2}}$ does exist. The result then follows by applying the same coupling scheme as in the proof of Lemma 2.1; see Neumann and Paparoditis (2005) for details. $\square$

**Proof of Theorem 4.1.** The method of proof is exactly the same as that for Theorem 3.1. Lemma 4.3 ensures that $(X_t^*)_{t\in\mathbb{Z}}$ satisfies appropriate conditions of weak dependence which yields, in conjunction with the fact that $P^{X_t^*}(B)$ converges to $P^{X_t}(B)$ with a sufficiently fast rate, that $(U_n^*)_{t\in\mathbb{Z}}$ is stochastically equicontinuous. Convergence of the finite-dimensional distributions to a Gaussian limit again follows from Theorem 6.1, while the result of Corollary 4.1 ensures that its covariance function is the same as for the original process under the null; for details, see Neumann and Paparoditis (2005). $\square$

## Acknowledgements

We would like to thank the Editor, the Associate Editor and two referees for their valuable comments leading to an improvement of the paper. We also thank Miguel A. Delgado for his helpful comments.